\documentclass[11pt,a4paper,reqno]{amsart}
\usepackage[T1]{fontenc}
\usepackage{lmodern}
\usepackage{amssymb}
\usepackage{geometry,booktabs,array,microtype,caption}
\usepackage{needspace}
\usepackage[numbers,sort&compress]{natbib}
\usepackage[hyperfootnotes=false, colorlinks, linkcolor={blue}, citecolor={magenta}, filecolor={blue}, urlcolor={blue}, plainpages=false, pdfpagelabels]{hyperref}
\allowdisplaybreaks[2]
\newtheorem{theorem}{Theorem}[section]
\newtheorem{proposition}[theorem]{Proposition}
\newtheorem{lemma}[theorem]{Lemma}
\newtheorem{corollary}[theorem]{Corollary}
\theoremstyle{remark}
\newtheorem{remark}[theorem]{Remark}
\title{On the moduli of distinct covering systems}
\author{Xiangjun Dai}
\address{School of Mathematical Sciences and LPMC, Nankai University,
Tianjin 300071, China}
\email{2120230022@mail.nankai.edu.cn}

\author{Huixi Li}
\address{School of Mathematical Sciences and LPMC, Nankai University,
Tianjin 300071, China}
\email{lihuixi@nankai.edu.cn}

\author{Junfeng Li}
\address{School of Mathematical Sciences and LPMC, Nankai University,
Tianjin 300071, China}
\email{junfengli.math@gmail.com}

\author{Shaoyun Yi}
\address{School of Mathematical Sciences, Xiamen University,
Xiamen, Fujian 361005, China}
\email{yishaoyun926@xmu.edu.cn}
\date{}
\hypersetup{pdftitle={On the moduli of distinct covering systems}}

\newcommand{\SquarefreeFiniteUpper}{0.999106529736774014}
\newcommand{\SquarefreeThirdMomentUpper}{145782367956.368914}
\newcommand{\SquarefreeTailUpper}{0.000047994195211974}
\newcommand{\SquarefreeTotalNumerator}{999154523931985986468684433520377470012137411}
\newcommand{\SquarefreeTotalUpper}{0.999154523931985987}
\newcommand{\OddFifteenEarlyUpper}{0.379554469493876041}
\newcommand{\OddFifteenFirstMomentUpper}{2.864960189566640836}
\newcommand{\OddFifteenThirdMomentUpper}{332.251231299406525656}
\newcommand{\OddFifteenFiniteUpper}{0.555492278736565973}
\newcommand{\OddFifteenMiddleUpper}{0.005321270446557438}
\newcommand{\OddFifteenTailUpper}{0.000079123607191286}
\newcommand{\OddFifteenTotalUpper}{0.940447142284190738}
\newcommand{\OddElevenEarlyUpper}{0.324201008245809909}
\newcommand{\OddElevenFiniteUpper}{0.663381303631749998}
\newcommand{\OddElevenMiddleUpper}{0.006095256428232006}
\newcommand{\OddElevenTailUpper}{0.000090632242845945}
\newcommand{\OddElevenTotalUpper}{0.993768200548637858}
\newcommand{\OddThirteenEarlyUpper}{0.388465921880343297}
\newcommand{\OddThirteenFiniteUpper}{0.567538808280385359}
\newcommand{\OddThirteenMiddleUpper}{0.005510192912603654}
\newcommand{\OddThirteenTailUpper}{0.000081932753455151}
\newcommand{\OddThirteenTotalUpper}{0.961596855826787461}
\newcommand{\NineFiniteUpper}{0.995566393703960831}
\newcommand{\NineMiddleUpper}{0.003831404539065249}
\newcommand{\NineTailUpper}{0.000045016154141302}
\newcommand{\NineTotalUpper}{0.999442814397167382}
\newcommand{\FiveFiniteUpper}{0.771326775551603255}
\newcommand{\FiveMiddleUpper}{0.065156698062518127}
\newcommand{\FiveTailUpper}{0.000008082310506662}
\newcommand{\FiveTotalUpper}{0.836491555924628044}

\begin{document}
\begin{abstract}
In 2022, Balister, Bollob\'as, Morris, Sahasrabudhe, and Tiba proved that the least common multiple of the moduli in any distinct covering system with odd moduli is divisible by $9$ or $15$. We strengthen this necessary condition to divisibility by $45$ and prove that the minimum modulus of an odd distinct covering system is at most $9$. For arbitrary distinct moduli, we lower the upper bound of Balister et al.\ on the minimum modulus from $616000$ to $14000$. For distinct squarefree moduli, we lower the upper bound of $118$ obtained by Cummings, Filaseta, and Trifonov in 2025 to $7$. Our proofs combine the distortion method with a joint comparison based on directional convexity. The comparison uses independent random exponents to bound expectations of convex functions of sums of indicators of residue classes, while retaining the restrictions on the allowed moduli. A sharper bound for the conditional density on each fibre further improves the estimates for the subsequent stages of the distortion process.
\end{abstract}
\maketitle

\section{Introduction and main results}\label{sec:intro}

A \emph{covering system} is a finite collection of residue classes
$a_i+d_i\mathbb{Z}$, indexed by $1\le i\le k$, whose union is $\mathbb{Z}$,
where $a_i\in\mathbb{Z}$ and the $d_i$ are integers greater than $1$.
It is called \emph{distinct} if $d_1,\ldots,d_k$ are pairwise distinct.
In this case, let $\mathcal{D}=\{d_1,\ldots,d_k\}$.
For each $d\in\mathcal{D}$, denote the unique residue class
in the system with modulus $d$ by $a_d+d\mathbb{Z}$.
We write $m=\min\mathcal{D}$ for its minimum modulus and
$Q=\operatorname{lcm}(\mathcal{D})$ for the least common multiple
of its moduli.

Erd\H{o}s~\cite{Erdos1950} introduced covering systems in 1950 in connection
with integers of the form $2^k+p$ and posed the \emph{minimum modulus problem},
asking whether $m$ can be arbitrarily large.
A second basic question, posed by Erd\H{o}s and Selfridge, asks whether every distinct covering system must contain an even modulus. The \emph{Erd\H{o}s--Selfridge conjecture} asserts an affirmative answer (see~\cite{Schinzel1967,Erdos1995}). 
For surveys and historical background, see Balister~\cite{Balister2024} and the
Erd\H{o}s Problems website~\cite{BloomErdosProblems}.

Hough and Nielsen~\cite{HoughNielsen2019} proved in 2019 that every
distinct covering system contains a modulus divisible by $2$ or $3$,
so that $3\mid Q$ whenever all moduli are odd.
Balister, Bollob\'as, Morris, Sahasrabudhe, and Tiba~\cite[Theorem~1.4]{BBMST2022}
subsequently established the stronger necessary condition $9\mid Q$ or $15\mid Q$.
Our first theorem shows that both divisibility conditions hold simultaneously.

\begin{theorem}\label{thm:odd-lcm}
If every modulus in a distinct covering system is odd, then $45\mid Q$.
\end{theorem}

For odd squarefree moduli, successive works of Berger, Felzenbaum, and Fraenkel~\cite{BFF1986,BFF1987} (see also Sun~\cite{Sun1991}),
Simpson and Zeilberger~\cite{SimpsonZeilberger1991}, and Guo and Sun~\cite{GuoSun2005}
showed that $Q$ would need at least $11$, $18$, and $22$ distinct prime factors,
before Balister et al.~\cite[Theorem~1.1]{BBMST2021squarefree} proved in 2021
that no odd squarefree distinct covering system exists.
Theorem~\ref{thm:odd-lcm} immediately yields an alternative proof of their theorem since $9\mid Q$.

We also obtain an explicit upper bound on the minimum modulus of any odd distinct covering system.

\begin{theorem}\label{thm:odd}
If every modulus in a distinct covering system is odd, then at least
one modulus is at most $9$.
\end{theorem}

We remark that odd coverings do exist once a single odd modulus is permitted to repeat. Harrington, Sun, and Wong~\cite{HarringtonSunWong2022} constructed odd coverings in which $7$ occurs four times or a prime $p\ge23$ occurs $p-5$ times, and Bispels et al.~\cite{BispelsEtAl2026} showed that $9$ may occur three times and extended the prime range to $p\ge17$, with all remaining moduli distinct.

Without the parity restriction, the following table records successive constructions up to Owens's example with minimum modulus $42$.
\begin{center}
\begin{tabular}{@{}llll@{}}
\toprule
Least modulus & Author & Reference & Year \\
\midrule
$6$  & Swift      & \cite{Swift1954} & 1954 \\
$9$  & Churchhouse & \cite{Churchhouse1968} & 1968 \\
$18$ & Krukenberg & \cite{Krukenberg1971} & 1971 \\
$20$ & Choi       & \cite{Choi1971} & 1971 \\
$24$ & Morikawa   & \cite{Morikawa1981examples} & 1981 \\
$25$ & Gibson     & \cite{Gibson2009} & 2009 \\
$40$ & Nielsen    & \cite{Nielsen2009} & 2009 \\
$42$ & Owens      & \cite{Owens2014} & 2014 \\
\bottomrule
\end{tabular}
\end{center}

In the opposite direction, Filaseta, Ford, Konyagin, Pomerance, and
Yu~\cite{FFKPY2007} proved in 2007 that $m$ is bounded whenever
$\sum_{d\in\mathcal{D}}1/d$ is bounded.  Hough~\cite{Hough2015} made a breakthrough and resolved
Erd\H{o}s's minimum modulus problem in 2015 by proving that $m\le10^{16}$.
Simplifying and sharpening Hough's approach, Balister et al.~\cite{BBMST2020, BBMST2022} developed the \emph{distortion method} to estimate the density of the uncovered set. They used this method to lower the general upper bound to $m<616000$~\cite[Theorem~8.1]{BBMST2022}. In a further striking application, they proved Schinzel's conjecture~\cite{Schinzel1967} that one modulus divides another~\cite[Theorem~1.2]{BBMST2022}. For distinct squarefree moduli, Cummings, Filaseta, and Trifonov~\cite[Theorem~1.1]{CFT2025} in 2025 obtained $m\le118$. Our next theorem improves both bounds.

\begin{theorem}\label{thm:bounds}
Every distinct covering system satisfies the following bounds.
\begin{enumerate}
\renewcommand{\theenumi}{\thetheorem(\roman{enumi})}
\renewcommand{\labelenumi}{(\roman{enumi})}
\item\label{thm:general} The least modulus is less than $14000$.
\item\label{thm:squarefree} If all moduli are squarefree, then the least
modulus is at most $7$.
\end{enumerate}
\end{theorem}

\begin{remark}\label{rem:numerical-threshold}
The proof of $m<14000$ leaves a margin exceeding $0.003$ in the final inequality ruling out a covering.
A further computer-assisted calculation, which keeps more detailed information about the preceding prime divisors, gives the stronger bound $m\le13278$.
In the squarefree case, a separate computer-assisted refinement, which tracks intersections of residue classes at small primes and uses refined convex estimates for subsequent prime stages, indicates that the bound can be improved to $m\le6$.
\end{remark}

Several related developments concern other restrictions on the moduli and analogues over other rings. For integer coverings with ordered moduli $d_1<d_2<\cdots$, Cummings, Filaseta, and Trifonov~\cite[Theorem~1.2]{CFT2025} bounded $d_j$ whenever the first $j-1$ classes do not already cover $\mathbb Z$. Klein, Koukoulopoulos, and Lemieux~\cite{KKL2024} gave quantitative bounds on $d_j$ for minimal coverings and bounded the minimum modulus under bounded multiplicity, while Gyarmati~\cite{Gyarmati2026} gave an elementary maximality argument deducing the boundedness of $d_2$ from that of $d_1$. Further results concern the sum of the reciprocals of the moduli. Filaseta and Kalogirou~\cite{FK2026} proved that there is an absolute $\varepsilon>0$ such that every distinct covering system with minimum modulus at least $5$ satisfies $\sum_{d\in\mathcal D}d^{-1}\ge1+\varepsilon$. More recently, Filaseta, Kalogirou, and Klein~\cite{FKK2026} obtained the explicit bound $1+10^{-28}$ when the minimum modulus is a prime $q\ge11$.

In the setting of number fields and global function fields, Li, Wang, and Yi~\cite{LiWangYi2023pre} resolved a number-field analogue of the minimum modulus problem, while Li, Wang, Wang, and Yi~\cite{LWWY2024,LWWY2025} settled the problem over $\mathbb{F}_q[x]$, disproved Azlin's conjecture~\cite{Azlin2011thesis}, and ruled out distinct coverings of $\mathbb{F}_q[x]$ for $q\ge759$. This bound was subsequently lowered to $q>73$ by Wang~\cite{BiaoWang2026} and to $q\ge19$ by Dai~\cite{Dai2026}. Dai also bounded the minimum degree over $\mathbb{F}_2[x]$ by $25$. See also Nie~\cite{Nie2026} for explicit constructions of covering systems of $\mathbb{F}_q[x]$ and Wang~\cite{RongyinWang2026,RongyinWang2026asymptotic} for finite and asymptotic results on covering systems of $\mathbb{F}_q[x]$.

\medskip
\noindent\textbf{Notation.}
For a positive integer $n$, let $\tau(n)$ denote the number of positive
divisors of $n$, $\Omega(n)$ the number of prime factors of $n$ counted with multiplicity, and $v_p(n)$ the exponent of $p$ in the prime factorization of $n$.
For $a\ge0$ and $n\ge1$, put
\begin{equation}\label{eq:small-divisors-and-support}
   D_a(n)=\#\{d\ge1:d\mid n,\ d<a\},\qquad \operatorname{supp}(n)=\{p:p\text{ is prime and }p\mid n\}.
\end{equation}
For $x\in\mathbb R$, write $x_+=\max\{x,0\}$, and put
$\mathbb R_+=[0,\infty)$. The indicator of an event $E$ is
$\mathbf1_E$.
We write $\mathbb{P}$ and $\mathbb{E}$ for probability and expectation, indicating the underlying
measure by a subscript when needed. When stages are indexed by primes,
$p^-$ denotes the stage immediately before $p$ in the chosen order,
so $\mu_{p^-}$ is the measure before the $p$-step.

\medskip
\noindent\textbf{Proof strategy.} Our approach builds on Hough's filtration by prime factors~\cite{Hough2015} and the distortion method of Balister, Bollob\'as, Morris, Sahasrabudhe, and Tiba~\cite{BBMST2022}. The latter uses a reweighted probability measure to control the mass of the newly covered residue classes. We develop this framework in two directions. First, we construct a random model that takes into account which moduli are allowed and use it to bound the contribution at each prime. Second, we use a joint convex comparison to relate the contribution at the current prime to those at later stages. This allows us to combine these improvements across stages. Together with rigorous finite recursions and analytic tail estimates, these refinements lead to the explicit bounds established in this paper.

More precisely, the distortion method excludes a covering when the sum of the stage contributions is less than $1$. At a prime $p$, the contribution is
\[c_p\mathbb E(\alpha_p-\delta_p)_+,\]
where $\alpha_p$ is the covered fibre proportion, $0\le\delta_p\le1/2$, and $c_p=(1-\delta_p)^{-1}$. We bound this expectation by comparing the congruence indicators in an upper bound for $\alpha_p$ with divisibility indicators defined by independent random exponents. This comparison retains the restrictions on the allowed moduli before the positive part is taken. For the bound $m<14000$, we also use the smaller conditional density bound $(1-\alpha_p)^{-1}$ on fibres where $\alpha_p<\delta_p$. A joint convex comparison then allows the resulting improvements at later stages to be combined without double counting.

To prove Theorem~\ref{thm:odd-lcm} in Section~\ref{sec:odd-lcm}, we establish $9\mid Q$ and $5\mid Q$ separately. Assuming $9\nmid Q$, we adjoin a residue class modulo $3$ if necessary and restrict the covering to the other two residue classes modulo $3$. The two induced coverings share the classes whose moduli are coprime to $3$, while their remaining reduced moduli are disjoint. Applying the scalar convex comparison to each branch and averaging via Lemma~\ref{lem:averaging} shows that the average of their distortion sums is less than $1$. To prove $5\mid Q$, we condition on the complement of all classes with $3$-power moduli. By Lemma~\ref{lem:pure-three-conditioning}, this bounds the probability of each residue class modulo $3^e$ by $2/3^e$ and leaves only primes $p\ge7$, again yielding a total mass below $1$.

In Section~\ref{sec:odd}, we prove Theorem~\ref{thm:odd} by first excluding moduli below $15$ from the comparison sum to rule out $m\ge15$. The remaining cases $m=11$ and $m=13$ are then excluded by using the uniform reference measure on the complement of the class with modulus $q$, where $q\in\{11,13\}$, in the $q$-coordinate and applying Lemma~\ref{lem:reference-law}.

Theorem~\ref{thm:bounds} is established in Section~\ref{sec:bounds-proof}. In Section~\ref{sec:general}, assuming $m\ge14000$, we keep the coordinates corresponding to primes up to $23$ uniformly distributed and use Shearer's lemma to bound the union of the classes assigned to these coordinates. The choice of the initial uniform block balances its covering contribution against the subsequent growth of the comparison distribution and its moments. Enlarging the block does not necessarily improve the final bound. We then apply the joint comparison of Lemmas~\ref{lem:joint-comparison} and~\ref{lem:common-kernel-propagation} for $p\ge29$ to bring the total mass below $997/1000$. In Section~\ref{sec:squarefree}, which is independent of Section~\ref{sec:general}, squarefreeness allows the comparison to use independent Bernoulli variables. A finite recursion for their sum, combined with a third-moment tail bound, excludes $m>9$ and hence yields $m\le7$.

Appendix~\ref{app:computations} records the computational verification, fixed parameter tables, and rational exponential bounds.

\section{Distortion and scalar comparison}\label{sec:distortion-scalar-comparison}
In this section, we develop the tools used to estimate the stage
contributions. We begin with the distortion method and a scalar
comparison, then give the finite recursions and tail bounds needed
for the numerical estimates. We also record a form of Shearer's lemma
and a conditioning result for classes with moduli that are powers
of three.

\subsection{The distortion method}\label{sec:distortion}
Factor $Q$ as
\begin{equation*}
Q=\prod_{j=1}^s q_j^{b_j},
\end{equation*}
where the $q_j$ are distinct primes, listed in an arbitrary fixed order.
Under the isomorphism
\begin{equation*}
\mathbb{Z}/Q\mathbb{Z}\xrightarrow{\sim}\prod_{j=1}^s\mathbb{Z}/q_j^{b_j}\mathbb{Z},
\end{equation*}
a residue class modulo $d\mid Q$ restricts only those coordinates $j$ for which $q_j\mid d$.

Let $X_0$ consist of a single point and put
\begin{equation*}
X_j=\prod_{i=1}^j\mathbb{Z}/q_i^{b_i}\mathbb{Z}.
\end{equation*}
Assign each original congruence class to the last coordinate on which it depends, and let $B_j\subseteq X_j$ denote the union of all classes assigned to the $j$-th coordinate. For $k>j$, we view $B_j$ as a cylinder set in $X_k$, unrestricted in later coordinates.
For $x\in X_{j-1}$, write
\begin{equation*}
F_j(x)=\{(x,y):y\in\mathbb{Z}/q_j^{b_j}\mathbb{Z}\},
\qquad
\alpha_j(x)=\frac{|F_j(x)\cap B_j|}{q_j^{b_j}}\in[0,1]
\end{equation*}
for the fibre above $x$ and the covered proportion of that fibre.

Given a probability measure $\mu_{j-1}$ on $X_{j-1}$, extend it uniformly to $X_j$ by
\begin{equation*}
\widetilde\mu_{j-1}(x,y)=\frac{\mu_{j-1}(x)}{q_j^{b_j}}
\qquad
(x\in X_{j-1},\ y\in\mathbb{Z}/q_j^{b_j}\mathbb{Z}).
\end{equation*}
Choose a distortion parameter $0\le\delta_j\le1/2$ and set
\begin{equation*}
c_j=\frac1{1-\delta_j}.
\end{equation*}
Following Balister et al.~\cite[Section~2.2]{BBMST2022}, define $\mu_j$ on $X_j$ by
\begin{equation}\label{eq:distorted-measure}
\mu_j(x,y)=
\begin{cases}
\displaystyle
\max\!\left\{0,\frac{\alpha_j(x)-\delta_j}
{\alpha_j(x)(1-\delta_j)}\right\}
\cdot \widetilde\mu_{j-1}(x,y), &(x,y)\in B_j,\\[9pt]
\displaystyle
\min\!\left\{\frac1{1-\alpha_j(x)},\frac1{1-\delta_j}\right\}
\cdot \widetilde\mu_{j-1}(x,y), &(x,y)\notin B_j.
\end{cases}
\end{equation}
Each branch of \eqref{eq:distorted-measure} is evaluated only when the corresponding subset of $F_j(x)$ is nonempty, so the multipliers are well defined. In the squarefree case, the same construction can be expressed recursively in terms of fibre weights, as in Cummings, Filaseta, and Trifonov~\cite{CFT2025}.

\begin{lemma}\label{lem:distortion}
For every $x\in X_{j-1}$,
\begin{equation}\label{eq:fiber-marginal}
\sum_{y\in\mathbb{Z}/q_j^{b_j}\mathbb{Z}}\mu_j(x,y)=\mu_{j-1}(x).
\end{equation}
Moreover,
\begin{equation}\label{eq:stage-measure}
\mu_j(B_j)=c_j\mathbb{E}_{\mu_{j-1}}\left(\alpha_j-\delta_j\right)_+,
\end{equation}
and
\begin{equation}\label{eq:pointwise-distortion}
\mu_j(x,y)\le c_j\widetilde\mu_{j-1}(x,y).
\end{equation}
Consequently, conditional on $x$ with $\mu_{j-1}(x)>0$, a prescribed residue class modulo $q_j^e$, with $1\le e\le b_j$, in the new coordinate has conditional probability at most $c_jq_j^{-e}$.
\end{lemma}

\begin{proof}
Fix $x$ and write $\alpha=\alpha_j(x)$ and $\delta=\delta_j$. If $\alpha\le\delta$, the first line of \eqref{eq:distorted-measure} vanishes and the multiplier on the complement of $B_j$ is $(1-\alpha)^{-1}$. Hence the total mass of the fibre is
\begin{equation*}
(1-\alpha)\frac1{1-\alpha}\,\mu_{j-1}(x)=\mu_{j-1}(x).
\end{equation*}
If $\alpha>\delta$, the two parts of the fibre have total coefficient
\begin{equation*}
\alpha\frac{\alpha-\delta}{\alpha(1-\delta)}
+(1-\alpha)\frac1{1-\delta}=1,
\end{equation*}
which again gives \eqref{eq:fiber-marginal}. The $\mu_j$-mass of $F_j(x)\cap B_j$ is respectively $0$ and
\begin{equation*}
\frac{\alpha-\delta}{1-\delta}\,\mu_{j-1}(x),
\end{equation*}
and summing over $x$ proves \eqref{eq:stage-measure}.

The two multipliers in \eqref{eq:distorted-measure} are at most $c_j$, giving \eqref{eq:pointwise-distortion}. A residue class modulo $q_j^e$ occupies exactly the proportion $q_j^{-e}$ of the fibre. Summing \eqref{eq:pointwise-distortion} over that residue class gives the final assertion.
\end{proof}

For simplicity, write $L_j:=\mu_j(B_j)$ for the contribution at the $j$-th stage. When stages are indexed directly by primes $p=q_j$, we replace the subscript $j$ by $p$ and write $B_p,\alpha_p,\delta_p,c_p,\mu_p$, and $L_p:=\mu_p(B_p)$.

\begin{corollary}\label{cor:criterion}
Let $\mu=\mu_s$. For every $j$ one has $\mu(B_j)=\mu_j(B_j)$. Consequently, if
\begin{equation}\label{eq:union-criterion}
\sum_{j=1}^s\mu_j(B_j)<1,
\end{equation}
then the original family of residue classes does not form a covering system.
\end{corollary}

\begin{proof}
Every subsequent distortion step preserves the joint distribution of the preceding coordinates. Since $B_j$ depends only on the first $j$ coordinates, its mass never changes afterward. Therefore
\begin{equation*}
\mu\left(\bigcup_jB_j\right)\le\sum_j\mu_j(B_j)<1,
\end{equation*}
so at least one residue class modulo $Q$ remains uncovered.
\end{proof}

The criterion in \eqref{eq:union-criterion} also has a grouped form. If $A=\bigcup_{j=1}^h B_j$ is the union of all original classes depending only on the first $h$ coordinates, then
\begin{equation*}
\mu(A)+\sum_{j>h}\mu_j(B_j)<1
\end{equation*}
likewise rules out a covering.

Since subsequent steps preserve the joint distribution of the first $j$ coordinates, the conditional probability bound in Lemma~\ref{lem:distortion} continues to hold under the final measure $\mu$, and we may combine it across all prime coordinates dividing a modulus.

\begin{lemma}\label{lem:residue-class-bound}
Suppose that the first $h$ coordinates have been processed and $d\mid\prod_{j\le h}q_j^{b_j}$. Then for every integer $a$,
\begin{equation}\label{eq:residue-class-bound}
\mu(a+d\mathbb{Z})\le\frac1d\prod_{q_j\mid d}c_j.
\end{equation}
\end{lemma}

\begin{proof}
Apply the conditional probability bound first in the last coordinate involved in $d$, and then successively to the preceding coordinates in reverse order. Each distinct prime contributes one factor $c_j$, whereas its complete exponent in $d$ contributes $q_j^{-v_{q_j}(d)}$. Multiplication gives \eqref{eq:residue-class-bound}.
\end{proof}

\begin{lemma}\label{lem:reference-law}
Let $\nu_q$ be a probability measure on the new prime coordinate,
and use the reference measure $\mu_{j-1}\otimes\nu_q$.
Let $\alpha_q$ be the $\nu_q$-measure of the covered part of each fibre.
If every class modulo $q^e$ has $\nu_q$-measure at most
$\eta_q q^{-e}$, then its conditional probability after distortion
is at most $\widehat c_q q^{-e}$, where
\begin{equation*}
\widehat c_q=\frac{\eta_q}{1-\delta_q}.
\end{equation*}
The preceding marginal is preserved, and \eqref{eq:stage-measure}
holds with $c_q=(1-\delta_q)^{-1}$.
\end{lemma}
\begin{proof}
Use the reference measure $\widetilde\mu_{j-1}(x,y)=\mu_{j-1}(x)\nu_q(y)$.
Define $\mu_j$ by \eqref{eq:distorted-measure} on its support,
and set $\mu_j=0$ elsewhere.
For a preceding point $x$, put
\[
B_x=\{y:(x,y)\in B_j\},\qquad
\alpha=\nu_q(B_x)=\alpha_q(x),\qquad
\delta=\delta_q,\qquad c=c_q.
\]
Summing the two parts of the fibre gives
\[
\begin{aligned}
\sum_y\mu_j(x,y)
&=\mu_{j-1}(x) \cdot
\begin{cases}
\displaystyle(1-\alpha)\frac{1}{1-\alpha},
    &\alpha\le\delta,\\[6pt]
c(\alpha-\delta)+c(1-\alpha),
    &\alpha>\delta
\end{cases}\\
&=\mu_{j-1}(x).
\end{aligned}
\]
Thus the preceding marginal is preserved.
The covered part has mass
\[
\sum_{y\in B_x}\mu_j(x,y)
=c\,\mu_{j-1}(x)(\alpha-\delta)_+.
\]
Consequently,
\[
\mu_j(B_j)
=c_q\sum_x\mu_{j-1}(x)
       (\alpha_q(x)-\delta_q)_+
=c_q\mathbb E_{\mu_{j-1}}
       (\alpha_q-\delta_q)_+,
\]
which is \eqref{eq:stage-measure}.
Both multipliers in \eqref{eq:distorted-measure} are at most
$c_q$, so
\[
0\le\mu_j(x,y)\le c_q\mu_{j-1}(x)\nu_q(y).
\]
For a prescribed class $C$ modulo $q^e$ and
$\mu_{j-1}(x)>0$, it follows that
\[
\mathbb P_{\mu_j}(y\in C\mid x)
=\frac{\sum_{y\in C}\mu_j(x,y)}{\mu_{j-1}(x)}
\le c_q\nu_q(C)
\le c_q\eta_q q^{-e}
=\widehat c_q q^{-e}.
\]
Any covering of the full residue set also covers the support of the reference measure, so Corollary~\ref{cor:criterion} applies.
\end{proof}

\subsection{Convex comparison and moment estimates}\label{sec:convex}

For fixed marginal probabilities, the expectation of an increasing
convex function of a nonnegative weighted sum of indicators is
maximized when the events are nested.
This is the Bernoulli case of the Lorentz--Tchen rearrangement
inequality~\cite[Section~4.3, Theorem~5]{Tchen1980}.

\begin{lemma}\label{lem:nested-indicator-comparison}
Let $I_1,\ldots,I_r$ be $\{0,1\}$-valued random variables with $\mathbb{P}(I_j=1)\le u_j\le1$. Let $U$ be uniform on $[0,1]$ and put $I_j^*=\mathbf{1}_{\{U\le u_j\}}$. If $b_j\ge0$ and $\Phi$ is increasing and convex, then
\begin{equation}\label{eq:nested-indicator-comparison}
\mathbb{E}\Phi\left(\sum_{j=1}^rb_jI_j\right)
\le
\mathbb{E}\Phi\left(\sum_{j=1}^rb_jI_j^*\right).
\end{equation}
\end{lemma}

\begin{proof}
Relabel the indices so that
\[
u_1\ge u_2\ge\cdots\ge u_r,
\qquad u_{r+1}=0.
\]
Since $\Phi$ is increasing and convex, we have
\begin{equation}\label{eq:convex-increments}
  \Phi(x+h)-\Phi(x)
\le \Phi(y+h)-\Phi(y)
\qquad (0\le x\le y,\ h\ge0).
\end{equation}
Define the partial sums
\[
s_0=Z_0=0,\qquad
s_j=\sum_{k=1}^j b_k,\qquad
Z_j=\sum_{k=1}^j b_k I_k
\quad (1\le j\le r),
\]
and put
\[
\Delta_j=\Phi(s_j)-\Phi(s_{j-1})\ge0.
\]
Since $0\le Z_{j-1}\le s_{j-1}$ and $I_j\in\{0,1\}$,
\eqref{eq:convex-increments} yields
\[
\Phi(Z_j)-\Phi(Z_{j-1})
=I_j\bigl\{\Phi(Z_{j-1}+b_j)-\Phi(Z_{j-1})\bigr\}
\le I_j\bigl\{\Phi(s_{j-1}+b_j)-\Phi(s_{j-1})\bigr\}
=I_j\Delta_j.
\]
Summing over $j$ and taking expectations, we obtain
\[
\mathbb E\Phi\left(\sum_{j=1}^r b_jI_j\right)
\le \Phi(0)+\sum_{j=1}^r \Delta_j\,\mathbb P(I_j=1)
\le \Phi(0)+\sum_{j=1}^r u_j\Delta_j.
\]
For the nested indicators $I_j^*=\mathbf1_{\{U\le u_j\}}$,
exactly the first $j$ indicators equal $1$ when
$u_{j+1}<U\le u_j$, while all indicators vanish when $U>u_1$.
Therefore
\[
\begin{aligned}
\mathbb E\Phi\left(\sum_{j=1}^r b_jI_j^*\right)
&=(1-u_1)\Phi(0)
  +\sum_{j=1}^r (u_j-u_{j+1})\Phi(s_j)\\
&=\Phi(0)+\sum_{j=1}^r
  u_j\bigl\{\Phi(s_j)-\Phi(s_{j-1})\bigr\}\\
&=\Phi(0)+\sum_{j=1}^r u_j\Delta_j.
\end{aligned}
\]
Combining the two estimates proves \eqref{eq:nested-indicator-comparison}.
\end{proof}

Let $\mathcal P$ be a finite set of primes with a fixed order, and put
\[
X=\prod_{q\in\mathcal P}\mathbb Z/q^{a_q}\mathbb Z,
\qquad a_q\in\mathbb Z_{\ge1}.
\]
Let $\mu$ be a probability measure on $X$.
Assume that, conditional on any values of the preceding coordinates
having positive $\mu$-probability, each residue class modulo $q^e$
in the $q$-coordinate has probability at most $\gamma_q/q^e$,
where $1\le e\le a_q$ and $1\le\gamma_q\le q$.
Let the independent nonnegative integer-valued variables $T_q$ satisfy
\begin{equation}\label{eq:comparison-exponent-law}
\mathbb P(T_q\ge e)=\frac{\gamma_q}{q^e}\qquad(e\ge1),
\end{equation}
and put $R=\prod_{q\in\mathcal P}q^{T_q}$.

\begin{theorem}\label{thm:scalar-comparison}
For $j=1,\ldots,r$, let
\[
\begin{aligned}
C_j&=\left\{x\in X:
x_q\equiv a_{j,q}\pmod{q^{e_{j,q}}}
\text{ for every }q\in\mathcal P\right\},\\
m_j&=\prod_{q\in\mathcal P}q^{e_{j,q}},
\end{aligned}
\]
where $a_{j,q}\in\mathbb Z$ and
$e_{j,q}\in\{0,\ldots,a_q\}$.
An exponent $e_{j,q}=0$ imposes no restriction.
For $b_j\ge0$ and every increasing convex function $\Phi$,
\begin{equation}\label{eq:scalar-comparison}
\mathbb E_\mu\Phi\left(
\sum_{j=1}^r b_j\mathbf1_{\{x\in C_j\}}\right)
\le
\mathbb E\Phi\left(
\sum_{j=1}^r b_j\mathbf1_{\{m_j\mid R\}}\right).
\end{equation}
\end{theorem}

\begin{proof}
Write the primes in their prescribed order as $q_1,\ldots,q_s$.
Take $x$ with distribution $\mu$, independently of all the $T_q$,
and define
\[
I_{j,q}
=\mathbf1_{\{x_q\equiv a_{j,q}\bmod q^{e_{j,q}}\}},
\qquad
I_{j,q}^*=\mathbf1_{\{T_q\ge e_{j,q}\}}.
\]
For $0\le k\le s$, put
\[
S_k=\sum_{j=1}^r b_j
       \prod_{i\le k}I_{j,q_i}
       \prod_{i>k}I_{j,q_i}^*,
\]
with empty products equal to $1$.

Fix $1\le k\le s$ and condition on
\[
x_{<k}=(x_{q_i})_{1 \le i < k},
\qquad T_{>k} = (T_{q_i})_{k < i \le s}. 
\]
Then
\[
\begin{gathered}
\beta_{j,k}
=b_j\prod_{i<k}I_{j,q_i}\prod_{i>k}I_{j,q_i}^*\ge0,\\
S_k=\sum_j\beta_{j,k}I_{j,q_k},
\qquad
S_{k-1}=\sum_j\beta_{j,k}I_{j,q_k}^*,
\end{gathered}
\]
where the weights $\beta_{j,k}$ are fixed under this conditioning.
Since $T_{>k}$ is independent of $x$, the assumed bound gives
\[
\mathbb P(I_{j,q_k}=1\mid x_{<k},T_{>k})
\le
\mathbb P(I_{j,q_k}^*=1\mid x_{<k},T_{>k})
=
\begin{cases}
1,&e_{j,q_k}=0,\\
\gamma_{q_k}q_k^{-e_{j,q_k}},&e_{j,q_k}>0.
\end{cases}
\]
The comparison indicators are nested, since they use the same
random exponent $T_{q_k}$.
Lemma~\ref{lem:nested-indicator-comparison}, applied with weights $\beta_{j,k}$,
therefore yields
\[
\mathbb E[\Phi(S_k)\mid x_{<k},T_{>k}]
\le
\mathbb E[\Phi(S_{k-1})\mid x_{<k},T_{>k}].
\]
Taking expectations and applying this inequality for
$k=s,s-1,\ldots,1$, we obtain
\[
\mathbb E_\mu\Phi\left(
\sum_j b_j\mathbf1_{\{x\in C_j\}}\right)
=\mathbb E\Phi(S_s)
\le\mathbb E\Phi(S_0)
=\mathbb E\Phi\left(
\sum_j b_j\mathbf1_{\{m_j\mid R\}}\right),
\]
which proves \eqref{eq:scalar-comparison}.
\end{proof}

We retain the unbounded random exponents $T_q$ in \eqref{eq:comparison-exponent-law},
even when the exponents in the congruence moduli are bounded.
Adding these divisors can only increase the upper bounds. If the reference coordinates are uniform, then $\gamma_q=c_q\le2\le q$.
The conditioned coordinates used below also satisfy $\gamma_q\le q$.

At stage $p$, let $m$ range over positive integers supported on
the preceding primes, and let $e\ge1$.
Under a uniform reference measure at $p$, a class modulo $mp^e$
occupies proportion $p^{-e}$ of the fibre when its congruence
modulo $m$ is satisfied.
The union bound therefore gives
\begin{equation*}
\alpha_p(x)\le
\sum_{mp^e\in\mathcal D}
p^{-e}\mathbf1_{\{x\equiv a_{mp^e}\bmod m\}}.
\end{equation*}
Here $R$ is formed from the preceding prime coordinates.
By distinctness, each pair $(m,e)$ occurs at most once, so
\begin{equation}\label{eq:divisor-upper}
\sum_{mp^e\in\mathcal D}
p^{-e}\mathbf1_{\{m\mid R\}}
\le
\sum_{m\mid R}\sum_{e\ge1}p^{-e}
=\frac{\tau(R)}{p-1}.
\end{equation}
Applying Theorem~\ref{thm:scalar-comparison} with
$\Phi(t)=(t-\delta_p)_+$ and using
\eqref{eq:stage-measure}, we obtain
\[
L_p=\mu_p(B_p)
=c_p\mathbb E_{\mu_{p^-}}(\alpha_p-\delta_p)_+
\le c_p\mathbb E\left(
\sum_{mp^e\in\mathcal D}
p^{-e}\mathbf1_{\{m\mid R\}}-\delta_p\right)_+
\le c_p\mathbb E\left(
\frac{\tau(R)}{p-1}-\delta_p\right)_+.
\]

For a nonuniform reference measure $\nu_p$ as in Lemma~\ref{lem:reference-law}, define $\alpha_p(x)$ as the $\nu_p$-measure of the covered part of the fibre above $x$, and replace each weight $p^{-e}$ above by $\eta_p p^{-e}$.
The stage factor remains $c_p=(1-\delta_p)^{-1}$.

\begin{lemma}\label{lem:averaging}
Let $A,X_1,\ldots,X_r,H\ge0$ satisfy
\begin{equation*}
A\le H,\qquad\sum_{j=1}^rX_j\le H.
\end{equation*}
Let $\lambda_j\ge0$, $\sum_j\lambda_j=1$, and write $\lambda_*=\max_j\lambda_j$. For every increasing convex function $\Phi$,
\begin{equation}\label{eq:averaging}
\sum_{j=1}^r\lambda_j\Phi(A+X_j)
\le(1-\lambda_*)\Phi(H)+\lambda_*\Phi(2H).
\end{equation}
\end{lemma}

\begin{proof}
By monotonicity we may first increase $A$ to $H$. With $H$ fixed, the left-hand side is a convex function of $(X_1,\ldots,X_r)$. The feasible vectors form the simplex whose vertices are the zero vector and the vectors $H e_j$. Hence convexity bounds the function by its maximum at a vertex.

At the nonzero vertex $H e_j$ the value is
\begin{equation*}
(1-\lambda_j)\Phi(H)+\lambda_j\Phi(2H).
\end{equation*}
Because $\Phi(2H)\ge\Phi(H)$, the largest such value corresponds to $\lambda_*$. The zero vertex is bounded by the same expression, proving \eqref{eq:averaging}.
\end{proof}

From \eqref{eq:comparison-exponent-law}, if $D_q=T_q+1$, then
\begin{equation*}
\mathbb{P}(D_q=1)=1-\frac{\gamma_q}{q},\qquad
\mathbb{P}(D_q=d)=\frac{\gamma_q(q-1)}{q^d}\quad(d\ge2).
\end{equation*}
For $r=1,2,3$, discrete summation gives
\begin{equation*}
\mathbb{E} D_q^r=1+\sum_{e\ge1}\bigl((e+1)^r-e^r\bigr)\mathbb{P}(T_q\ge e),
\end{equation*}
and hence
\begin{equation}\label{eq:moments}
\mathbb{E} D_q=1+\frac{\gamma_q}{q-1},\qquad
\mathbb{E} D_q^2=1+\gamma_q\frac{3q-1}{(q-1)^2},\qquad
\mathbb{E} D_q^3=1+\gamma_q\frac{7q^2-2q+1}{(q-1)^3}.
\end{equation}
Thus the untruncated moments can be computed exactly by recursion, independently of the truncated distribution table.

For every nonnegative random variable $X$ and $a\ge0$,
\begin{equation}\label{eq:positive-part}
\mathbb{E}(X-a)_+=\mathbb{E} X-a+\mathbb{E}(a-X)_+.
\end{equation}
Consequently, the positive-part expectation can be computed from the untruncated first moment together with the distribution on $X<a$.

For $x\ge0$ and $a>0$,
\begin{equation}\label{eq:cubic}
(x-a)_+\le\frac{4x^3}{27a^2}.
\end{equation}
The inequality is immediate for $x<a$. For $x\ge a$ it follows from
\begin{equation*}
4x^3-27a^2x+27a^3=(2x-3a)^2(x+3a)\ge0.
\end{equation*}

\subsection{Finite recursions and prime tail bounds}\label{sec:finite}
We first describe the rounding rules and the finite recursion used to estimate products of independent random variables. We then introduce a truncation that keeps the recursion finite.

Fix an integer scale $W\ge1$. For positive integers $r,s$,
\[
0\le\frac\ell W\le x\le\frac uW
\quad\Longrightarrow\quad
\frac1W\left\lfloor\frac{r\ell}{s}\right\rfloor
\le\frac rs x\le
\frac1W\left\lceil\frac{ru}{s}\right\rceil.
\]
For two-sided bounds $A^-\le A\le A^+$ and
$B^-\le B\le B^+$, we use
\[
A^{-}-B^{+}\le A-B\le A^{+}-B^{-}.
\]

Let $X,D$ be independent positive integer-valued variables,
and suppose
\[
W\mathbb P(X=t)\le P_t,\qquad
W\mathbb P(D=d)\le Q_d,
\]
where $P_t,Q_d$ are nonnegative integers.
Set
\begin{equation}\label{eq:product-recursion}
P'_u=\left\lceil
\frac1W\sum_{td=u}P_tQ_d
\right\rceil.
\end{equation}
Independence gives
\[
\mathbb P(XD=u)
=\sum_{td=u}\mathbb P(X=t)\mathbb P(D=d)
\le\frac1{W^2}\sum_{td=u}P_tQ_d
\le\frac{P'_u}{W}.
\]
Rounding each summand upward also gives an upper bound, since
\[
P'_u\le\sum_{td=u}\left\lceil\frac{P_tQ_d}{W}\right\rceil.
\]
Lower bounds use downward rounding, and the truncated probability distributions are never renormalized. For an integer cutoff $H\ge1$,
\[
td\le H,\quad t,d\ge1\quad\Longrightarrow\quad t\le H.
\]
Thus the multiplicative update on $u\le H$ uses only values at most $H$.

For a finite product $X=\prod_{i=1}^nD_i$ of independent
positive integer-valued multipliers, choose an integer $d_*\ge1$
and put
\[
\widetilde D_i=
\begin{cases}
D_i,&D_i\le d_*,\\
1,&D_i>d_*,
\end{cases}
\qquad
\widetilde X=\prod_{i=1}^n\widetilde D_i\le X.
\]
Since $(a-X)_+\le(a-\widetilde X)_+$,
\eqref{eq:positive-part} gives, for every real $a\ge0$,
\begin{equation}\label{eq:replacement-bound}
\begin{aligned}
\mathbb E(X-a)_+
&\le\mathbb EX-a+\mathbb E(a-\widetilde X)_+\\
&=\mathbb EX-a+
\sum_{1\le t<a}(a-t)\mathbb P(\widetilde X=t).
\end{aligned}
\end{equation}
When $a\le H$, the last sum uses only retained probabilities. For $k=1,2,3$, the full moments $\mathbb EX^k=\prod_{i=1}^n\mathbb ED_i^k$ are updated separately from the untruncated multipliers via~\eqref{eq:moments}.

The following lemma gives a finite recursion for a truncated product of random variables and uses it to bound the positive-part expectation of the original product.
\begin{lemma}\label{lem:whole-reset}
Let $X_0$ and $D_1,D_2,\ldots$ be independent positive integer-valued
random variables, and let $H,d_*\ge1$ be integers. Put $Y_0=X_0$ and
\begin{equation*}
X_j=X_{j-1}D_j,\qquad
Y_j=\begin{cases}
Y_{j-1}D_j,&D_j\le d_*,\\
1,&D_j>d_*.
\end{cases}
\end{equation*}
Then $1\le Y_j\le X_j$ almost surely. Writing
$\pi_j(d)=\mathbb{P}(D_j=d)$, $\rho_j=\mathbb{P}(D_j>d_*)$, and
$q_j(u)=\mathbb{P}(Y_j=u)$, the retained probabilities satisfy, for $1\le u\le H$,
\begin{equation}\label{eq:whole-reset-law}
q_j(u)=\sum_{\substack{1\le d\le d_*\\d\mid u}}
\pi_j(d)q_{j-1}(u/d)+\mathbf{1}_{\{u=1\}}\rho_j.
\end{equation}
If $\mathbb{E} X_j<\infty$ and $0\le a\le H$, then
\begin{equation}\label{eq:whole-reset-positive}
\mathbb{E}(X_j-a)_+\le\mathbb{E} X_j-a+
\sum_{\substack{1\le u\le H\\u<a}}(a-u)q_j(u).
\end{equation}
\end{lemma}

\begin{proof}
Starting from $Y_0=X_0\ge1$, induction gives
\[
\begin{array}{ll}
D_j\le d_*:&1\le Y_j=Y_{j-1}D_j\le X_{j-1}D_j=X_j,\\[2pt]
D_j>d_*:&1=Y_j\le X_j.
\end{array}
\]
Since $Y_{j-1}$ depends only on $X_0,D_1,\ldots,D_{j-1}$,
it is independent of $D_j$. Thus
\[
\mathbb P(Y_j=u,\ D_j\le d_*)
=\sum_{\substack{1\le d\le d_*\\d\mid u}}
\pi_j(d)q_{j-1}(u/d),
\]
whereas
\[
\mathbb P(Y_j=u,\ D_j>d_*)
=\mathbf1_{\{u=1\}}\cdot\rho_j\sum_{v\ge1}q_{j-1}(v)
=\mathbf1_{\{u=1\}}\cdot\rho_j.
\]
Adding proves \eqref{eq:whole-reset-law}. For $u \le H$, every divisor quotient $u/d$ is at most $H$, while the tail term $\mathbf{1}_{\{u=1\}}\rho_j$ accounts for all values of $Y_{j-1}$, including $Y_{j-1} > H$.

Finally, \eqref{eq:positive-part} and $Y_j\le X_j$ give
\[
\begin{aligned}
\mathbb E(X_j-a)_+
&=\mathbb EX_j-a+\mathbb E(a-X_j)_+\\
&\le\mathbb EX_j-a+\mathbb E(a-Y_j)_+\\
&=\mathbb EX_j-a+
\sum_{\substack{1\le u\le H\\u<a}}(a-u)q_j(u),
\end{aligned}
\]
since $a\le H$. This is \eqref{eq:whole-reset-positive}.
\end{proof}

We now give an integer version of this recursion for rigorous
computation. Let $W$ be a positive integer. Suppose that
$Q_{j-1}(u)$, $\Pi_j(d)$, and $P_j$ are nonnegative integer upper
bounds for $Wq_{j-1}(u)$, $W\pi_j(d)$, and $W\rho_j$, respectively.
Initially, choose $Q_0(u)\ge W\cdot\mathbb P(X_0=u)$ for $1\le u\le H$.
By induction, the update
\begin{equation}\label{eq:whole-reset-integer}
Q_j(u)=\sum_{\substack{1\le d\le d_*\\d\mid u}}
\left\lceil\frac{\Pi_j(d)Q_{j-1}(u/d)}{W}\right\rceil
+\mathbf{1}_{\{u=1\}}P_j
\end{equation}
satisfies $Q_j(u)/W\ge q_j(u)$ for $1\le u\le H$.
The expectation
\[
\mathbb E X_j=\mathbb E X_0\prod_{i=1}^j\mathbb E D_i
\]
is computed separately using the untruncated multipliers.
In the numerical verification, an upper bound for this expectation
is used in \eqref{eq:whole-reset-positive}.

The following lemma bounds the contribution from the remaining prime stages.
\begin{lemma}\label{lem:index-tail}
Let $N\ge1$ be an integer, and let $a\ge0$ and $c>0$ satisfy $a<2c$. Suppose that the nonnegative sequence $R_j$ satisfies
\begin{equation*}
R_j\le R_{N+1}\prod_{N<i<j}\left(1+\frac a{ci}\right)
\end{equation*}
for every $j > N$.
Let $C\ge0$. If the contribution at the $j$-th step is at most
$CR_j/(cj)^3$, then
\begin{equation}\label{eq:index-tail}
\sum_{j>N}\frac{CR_j}{(cj)^3}
\le\frac{CR_{N+1}}{c^2(2c-a)N^2}.
\end{equation}
\end{lemma}

\begin{proof}
Using $1+u\le e^u$ and the integral bound for the harmonic sum,
\begin{equation*}
\prod_{N<i<j}\left(1+\frac a{ci}\right)
\le\exp\left(\frac ac\sum_{N<i<j}\frac1i\right)
\le\left(\frac jN\right)^{a/c}.
\end{equation*}
Because $a/c<2$, the function $x^{a/c-3}$ is positive and decreasing. Hence
\begin{equation*}
\sum_{j>N}j^{a/c-3}
\le\int_N^\infty x^{a/c-3}\,dx
=\frac{N^{a/c-2}}{2-a/c}.
\end{equation*}
Substitution gives \eqref{eq:index-tail}.
\end{proof}

To apply Lemma~\ref{lem:index-tail} to prime stages, we first give an explicit lower bound for $p_j-1$ in terms of $j$.
Let $p_j$ denote the $j$-th prime. Dusart's explicit lower bound for primes \cite[Proposition~6.7]{Dusart2010} implies that, for every $j\ge664579$,
\begin{equation}\label{eq:prime-fifteen}
p_j-1>15j.
\end{equation}
The rational bounds
\begin{equation}\label{eq:exponential-bounds}
e^{13.406}<664579,\qquad e^{2.595}<13.406,
\end{equation}
verified in Appendix~\ref{app:exponential-checks}, give
$\log j>13.406$ and $\log\log j>2.595>2.1$. Hence
\[
\frac{p_j}{j}
\ge\log j+\log\log j-1+
\frac{\log\log j-2.1}{\log j}
>15.001>15+\frac1j.
\]
For a second tail estimate, Rosser's inequality $p_j>j\log j$~\cite{Rosser1939}
gives the following bound for $j\ge3$:
\begin{equation}\label{eq:rosser-shift}
p_j-1>p_{j-1}>(j-1)\log(j-1).
\end{equation}

The next lemma uses this bound to estimate the corresponding prime tail.
\begin{lemma}\label{lem:log-tail}
Suppose the computation has been carried through the $J$-th prime and
set $b=J-1$. Let $R_j\ge0$ be the third-moment bound before $p_j$,
with $R_{j+1}\le R_j(1+\beta/(p_j-1))$ for $j>J$ and $\beta\ge3$.
Let $C\ge0$ and let the stage contributions satisfy
$0\le L_j\le CR_j/(p_j-1)^3$.
If
\begin{equation*}
\log b\ge\ell>\frac{\beta-3}{2},
\end{equation*}
then
\begin{equation}\label{eq:log-tail}
\sum_{j>J}L_j
\le\frac{CR_{J+1}}{b^2\ell^2(2\ell-\beta+3)}.
\end{equation}
\end{lemma}

\begin{proof}
From \eqref{eq:rosser-shift} and an integral comparison,
\begin{equation*}
R_j\le R_{J+1}
\left(\frac{\log(j-1)}{\log b}\right)^\beta.
\end{equation*}
Since $\ell>(\beta-3)/2$, the function $(\log x)^{\beta-3}/x^3$ is decreasing on the relevant range. Hence the required sum is at most
\begin{equation*}
\frac{CR_{J+1}}{(\log b)^\beta}
\int_b^\infty\frac{(\log x)^{\beta-3}}{x^3}\,dx.
\end{equation*}
Put $x=be^t$ and use
\begin{equation*}
\left(1+\frac t{\log b}\right)^{\beta-3}
\le\exp\left(\frac{\beta-3}{\log b}t\right).
\end{equation*}
This yields \eqref{eq:log-tail} with $\log b$ in place of $\ell$. The denominator is increasing as a function of $\log b$ throughout the present range, so replacing $\log b$ by the smaller value $\ell$ preserves an upper bound.
\end{proof}

\subsection{Events with a dependency graph}\label{sec:dependency-graph}
We recall the following form of Shearer's lemma~\cite{Shearer1985} and include a proof for completeness.
\begin{lemma}\label{lem:shearer}
Let a finite family of events have dependency graph $G$, and suppose the probability of the event at vertex $v$ is at most $w_v$. For a vertex set $U$, define
\begin{equation*}
q(U)=\sum_{\substack{I\subseteq U\\I\text{ independent}}}
(-1)^{|I|}\prod_{v\in I}w_v.
\end{equation*}
If $q(U)>0$ for every $U$, then the probability that none of the events occurs is at least $q(V(G))$.
\end{lemma}

\begin{proof}
Let $F(U)$ be the probability that no event indexed by $U$ occurs. We prove by induction on $|U|$ that, for every $v\in U$,
\begin{equation}\label{eq:shearer-ratio}
\frac{F(U)}{F(U\setminus\{v\})}
\ge\frac{q(U)}{q(U\setminus\{v\})}>0.
\end{equation}
For a singleton this is $F(\{v\})\ge1-w_v$. Assume the assertion for smaller sets and let $N[v]$ be the closed neighborhood of $v$. The event at $v$ is jointly independent of all events at nonneighbors, so
\begin{equation*}
F(U)\ge F(U\setminus\{v\})-w_vF(U\setminus N[v]).
\end{equation*}
Multiplying the induction hypothesis along a chain from $U\setminus N[v]$ to $U\setminus\{v\}$ gives
\begin{equation*}
\frac{F(U\setminus N[v])}{F(U\setminus\{v\})}
\le\frac{q(U\setminus N[v])}{q(U\setminus\{v\})}.
\end{equation*}
Partitioning the independent sets according to whether they contain $v$ gives
\begin{equation*}
q(U)=q(U\setminus\{v\})-w_vq(U\setminus N[v]).
\end{equation*}
Combining these inequalities proves \eqref{eq:shearer-ratio}. Multiplying \eqref{eq:shearer-ratio} along a chain from the empty set to the full vertex set, and using $F(\varnothing)=q(\varnothing)=1$, gives the asserted inequality.
\end{proof}

\subsection{Classes with moduli that are powers of three}\label{sec:three-power-classes}
We first condition on the complement of the classes whose moduli are powers of three.
\begin{lemma}\label{lem:pure-three-conditioning}
Assume that an odd distinct covering system satisfies $5\nmid Q$ and put
$h=\max\{1,v_3(Q)\}$. After the classes whose moduli are powers of $3$
have been removed, the remaining classes cover a product probability space
in which every residue class modulo $3^e$, for $1\le e\le h$, has probability
at most $2/3^e$. Every other prime coordinate that can occur is indexed by a prime at least $7$.
\end{lemma}

\begin{proof}
Write $Q=3^{h_0}M$ with $(3,M)=1$, and use the period $Q^*=3^hM$.
If $h_0=0$, this only adds an unused uniform coordinate modulo $3$.
Every original modulus divides $Q^*$, and a covering modulo $Q$ also covers
modulo $Q^*$.
Let $P$ be the union of the original residue classes whose moduli are
exactly $3^e$. Distinctness gives at most one such class for each $e$.
For the uniform measure $u$ on $\mathbb{Z}/Q^*\mathbb{Z}$,
\begin{equation}\label{eq:pure-three-mass}
u(P)\le\sum_{e=1}^h\frac1{3^e}
=\frac12-\frac1{2\cdot3^h}<\frac12.
\end{equation}
Put $S=P^c$. The classes not used in the definition of $P$ cover $S$.

The set $S$ depends only on the coordinate $\mathbb{Z}/3^h\mathbb{Z}$.
Conditioning $u$ on $S$ therefore leaves a product of its conditioned
$3$-coordinate marginal and the uniform measure on $\mathbb{Z}/M\mathbb{Z}$.
For a residue class $C$ modulo $3^e$, where $1\le e\le h$,
equation \eqref{eq:pure-three-mass} gives
\begin{equation}\label{eq:conditioned-three-cap}
u(C\mid S)\le\frac{u(C)}{u(S)}\le\frac2{3^e}.
\end{equation}
The assumption $5\nmid Q$ and oddness show that every prime divisor of $M$ is at least $7$.
\end{proof}

The following corollary applies the distortion method after this conditioning.
\begin{corollary}\label{cor:conditioned-distortion}
In the setting of Lemma~\ref{lem:pure-three-conditioning},
let $\mu_0$ be the conditional probability measure on
$\mathbb Z/3^h\mathbb Z$ obtained in that lemma.
Starting with $\mu_0$, process the remaining prime coordinates
$p\ge7$ by uniform extension and distortion, with
$0\le\delta_p\le1/2$ and $c_p=(1-\delta_p)^{-1}$.
Then Theorem~\ref{thm:scalar-comparison} applies with
\[
\gamma_3=2,\qquad \gamma_p=c_p\quad(p\ge7).
\]
The remaining classes cannot cover $S$, the complement defined in Lemma~\ref{lem:pure-three-conditioning}, if
\[
\sum_{p\ge7}
c_p\mathbb E_{\mu_{p^-}}(\alpha_p-\delta_p)_+<1,
\]
where the sum runs over the remaining prime coordinates.
\end{corollary}

\begin{proof}
Let $S$ be the complement of the classes whose moduli are powers of $3$,
as in Lemma~\ref{lem:pure-three-conditioning}.
By \eqref{eq:conditioned-three-cap},
\[
\mu_0(a+3^e\mathbb Z)\le\frac2{3^e}
\qquad(1\le e\le h).
\]
At each subsequent prime $p$, Lemma~\ref{lem:distortion} gives
\[
\sum_y\mu_p(x,y)=\mu_{p^-}(x),\qquad
\mu_p(B_p)=c_p\mathbb E_{\mu_{p^-}}(\alpha_p-\delta_p)_+
\]
and
\[
\mathbb P_{\mu_p}(y\equiv a\bmod p^e\mid x)
\le\frac{c_p}{p^e}.
\]
The last bound holds for $1\le e\le v_p(Q)$ and
$\mu_{p^-}(x)>0$. Let $\mu$ be the final measure.
Preservation of the joint distribution of all preceding coordinates
keeps these conditional bounds valid under $\mu$. Since
\[
\gamma_3=2\le3,\qquad\gamma_p=c_p\le2\le p,
\]
Theorem~\ref{thm:scalar-comparison} applies.

The same marginal preservation gives
$\mu(B_p)=\mu_p(B_p)$ and $\mu(S)=1$. Consequently,
\[
\mu\left(S\setminus\bigcup_{p\ge7}B_p\right)
\ge1-\sum_{p\ge7}\mu(B_p)
=1-\sum_{p\ge7}
c_p\mathbb E_{\mu_{p^-}}(\alpha_p-\delta_p)_+.
\]
A positive right-hand side excludes a covering of $S$.
\end{proof}

\section{Joint comparison of stage contributions}\label{sec:joint-comparison}
The aim of this section is to bound the contributions of several stages jointly. We first establish a comparison for vectors of divisor weights and then use it to take advantage of a smaller conditional density bound.

\subsection{Comparison of divisor weight vectors}
Order the prime coordinates as in Section~\ref{sec:distortion}.
For $0\le\delta_p\le1/2$, put
\[
c_p=(1-\delta_p)^{-1},\qquad
g_p(a)=c_p(a-\delta_p)_+,\qquad
\mu_p(B_p)=\mathbb E_{\mu_{p^-}}g_p(\alpha_p).
\]
For a fixed base point $x$, \eqref{eq:distorted-measure} bounds
the conditional density multiplier by
\begin{equation*}
 \kappa_p(a)=
 \begin{cases}(1-a)^{-1},&0\le a<\delta_p,\\
 c_p,&a\ge\delta_p,
 \end{cases}
\end{equation*}
whenever $\alpha_p(x)\le a$. Thus any upper bound $\alpha_p\le a$ simultaneously yields the stage contribution bound $g_p(a)$ and the conditional density multiplier $\kappa_p(a)$ for subsequent coordinates.

Each modulus $mp^e\in\mathcal{D}$ contributes a nonnegative divisor weight $p^{-e}\mathbf{1}_{\{x\equiv a_{mp^e}\bmod m\}}$ determined by the preceding coordinates. By distinctness, each pair $(m,e)$ occurs at most once, and we may add the nonnegative terms corresponding to admissible moduli that are absent, obtaining an upper bound by summing over all admissible moduli.

To compare vector-valued divisor-weight sums across multiple primes, we work with \emph{increasing directionally convex} functions, namely functions $f:\mathbb{R}_+^N\to\mathbb{R}$ that are coordinatewise increasing and satisfy
\begin{equation}\label{eq:directional-convexity}
 f(x+u+v)-f(x+u)-f(x+v)+f(x)\ge0\qquad(u,v\ge0).
\end{equation}
For continuous functions, \eqref{eq:directional-convexity} is equivalent to separate
coordinatewise convexity combined with supermodularity.

\begin{lemma}\label{lem:vector-comparison}
In the setting of Theorem~\ref{thm:scalar-comparison}, let
$b_1,\ldots,b_r\in\mathbb R_+^N$ and let
$\Phi:\mathbb R_+^N\to\mathbb R$ be increasing and directionally
convex. Then
\[
\mathbb E_\mu\Phi\left(
\sum_{j=1}^rb_j\mathbf1_{\{x\in C_j\}}\right)
\le
\mathbb E\Phi\left(
\sum_{j=1}^rb_j\mathbf1_{\{m_j\mid R\}}\right).
\]
The same comparison integer $R$ is used in every component.
For coordinates obtained by distortion from a uniform reference measure, we take
$\gamma_p=c_p$, so that
\begin{equation*}
\mathbb P(T_p\ge e)=c_pp^{-e}\quad(e\ge1),
\qquad D_p=T_p+1.
\end{equation*}
\end{lemma}

\begin{proof}
By \eqref{eq:directional-convexity},
\[
\Phi(x+h)-\Phi(x)\le\Phi(y+h)-\Phi(y)
\qquad(0\le x\le y,\ h\ge0),
\]
where vector inequalities are componentwise.
First take indicators $I_j$ with $\mathbb P(I_j=1)\le u_j\le1$,
and relabel them so that $u_1\ge\cdots\ge u_r$.
Let $U$ be uniform on $[0,1]$ and set
\[
I_j^*=\mathbf1_{\{U\le u_j\}},\qquad
s_j=\sum_{\ell=1}^j b_\ell,\qquad
Z_j=\sum_{\ell=1}^j b_\ell I_\ell,\qquad s_0=Z_0=0.
\]
Since $0\le Z_{j-1}\le s_{j-1}$,
\[
\begin{aligned}
\Phi(Z_j)-\Phi(Z_{j-1})
&=I_j\{\Phi(Z_{j-1}+b_j)-\Phi(Z_{j-1})\}\\
&\le I_j\{\Phi(s_j)-\Phi(s_{j-1})\}.
\end{aligned}
\]
The increments on the right are nonnegative. The nested indicators
have an initial segment of active indices, so
\[
\Phi\left(\sum_j b_jI_j^*\right)
=\Phi(0)+\sum_j I_j^*\{\Phi(s_j)-\Phi(s_{j-1})\}.
\]
Consequently,
\[
\begin{aligned}
\mathbb E\Phi\left(\sum_j b_jI_j\right)
&\le\Phi(0)+\sum_j u_j\{\Phi(s_j)-\Phi(s_{j-1})\}\\
&=\mathbb E\Phi\left(\sum_j b_jI_j^*\right).
\end{aligned}
\]

Use the joint construction and indicators from the proof of
Theorem~\ref{thm:scalar-comparison}, with vector weights $b_j$, and put
\[
S_k=\sum_j b_j\prod_{i\le k}I_{j,q_i}
                       \prod_{i>k}I_{j,q_i}^*.
\]
Conditional on $x_{<k},T_{>k}$, the vectors
\[
\beta_{j,k}=b_j\prod_{i<k}I_{j,q_i}
                    \prod_{i>k}I_{j,q_i}^*\ge0
\]
are fixed. The conditional probability bounds are unchanged,
and the comparison indicators use the same $T_{q_k}$.
The one-step inequality therefore gives, almost surely,
\[
\mathbb E[\Phi(S_k)\mid x_{<k},T_{>k}]
\le\mathbb E[\Phi(S_{k-1})\mid x_{<k},T_{>k}].
\]
Taking expectations for $k=s,s-1,\ldots,1$ yields
\[
\mathbb E_\mu\Phi\left(\sum_jb_j\mathbf1_{\{x\in C_j\}}\right)
=\mathbb E\Phi(S_s)\le\cdots\le\mathbb E\Phi(S_0)
=\mathbb E\Phi\left(\sum_jb_j\mathbf1_{\{m_j\mid R\}}\right).
\]
\end{proof}
For later use, let $b\ge1$ be an integer. The distribution of
$T_p$ in Lemma~\ref{lem:vector-comparison} gives
\[
\mathbb P(T_p\wedge b=b)=\mathbb P(T_p\ge b)=c_pp^{-b}.
\]

\subsection{Preserving a smaller conditional density bound}

Fix a prime coordinate $p$, put $\delta=\delta_p$, and write
$c=(1-\delta)^{-1}$, $g(a)=c(a-\delta)_+$,
and set $\kappa=\kappa_p$ and $\Delta(a)=c-\kappa(a)$.
Let $\Phi$ be an increasing directionally convex functional that bounds the contribution of subsequent stages in terms of the remaining divisor weights, and let
$A_0,\ldots,A_b$ denote the nonnegative weight column vectors
indexed by the exponent of $p$ in the congruence moduli. Put
\begin{equation}\label{eq:weight-column-differences}
S_v=\sum_{e=0}^v A_e,\qquad
 B(A)=\sum_{v=1}^b p^{-v}\{\Phi(S_v)-\Phi(S_{v-1})\},
\end{equation}
and
\[
 V_\rho(A)=\Phi(A_0)+\rho B(A),\qquad1\le\rho\le c.
\]
Here $V_\rho(A)=\mathbb E_\rho\Phi(S_{T\wedge b})$ under
$\mathbb P_\rho(T\ge v)=\rho p^{-v}$ for $v\ge1$,
and $B(A)$ is the coefficient of $\rho$ in the subsequent-stage bound $V_\rho(A)$.

The following lemma combines the contribution of the current stage with an upper bound for the contributions of subsequent stages.
\begin{lemma}\label{lem:joint-comparison}
For any fixed $\vartheta$ with $0 \le \vartheta \le 1/c$, the function
\begin{equation}\label{eq:joint-bound}
  F(a, A) = g(a) + V_c(A) - \Delta(a)\min\{B(A), \vartheta\}
\end{equation}
is jointly increasing and directionally convex in $a$ and in all components of $A=(A_0,\ldots,A_b)$, and it bounds
the sum of the current stage contribution and the conditional contribution of subsequent stages from above.
\end{lemma}

\begin{proof}
First, $B(A)\ge0$ is coordinatewise increasing. Monotonicity in $A_0$ follows
from the increasing differences of $\Phi$, while for $A_j$ with $j\ge1$, summation by parts gives
\[
 B(A)=-p^{-1}\Phi(S_0)+\sum_{v=1}^{b-1}(p^{-v}-p^{-v-1})\Phi(S_v)
       +p^{-b}\Phi(S_b),
\]
in which every term depending on $A_j$ carries a positive coefficient.

Rewrite \eqref{eq:joint-bound} as the maximum of two expressions
\begin{equation}\label{eq:joint-bound-branches}
 F=\max\{g+V_\kappa,\;g+V_c-\Delta \vartheta\}.
\end{equation}
For fixed $a$, each branch is increasing, coordinatewise convex, and
supermodular in $A$, since $V_\rho$ is a convex combination of functions obtained by composing $\Phi$ with nonnegative linear maps. Write
\[
f_1=g+V_\kappa,\qquad
f_2=g+V_c-\Delta \vartheta,\qquad
\eta=f_2-f_1=\Delta(a)(B(A)-\vartheta).
\]
Since $\Delta(a)\ge0$ and $B$ is increasing, $\eta$ is coordinatewise increasing in $A$ for fixed $a$.

Put $w=x\wedge y$ and $z=x\vee y$.
If the same branch attains the maximum at both $x$ and $y$,
its supermodularity gives
$F(x)+F(y)\le F(w)+F(z)$.
Otherwise, after interchanging $x,y$ if necessary, suppose
$F(x)=f_2(x)$ and $F(y)=f_1(y)$.
Then
\[
\begin{aligned}
F(x)+F(y)
&=f_1(x)+f_1(y)+\eta(x)\\
&\le f_1(w)+f_1(z)+\eta(z)\\
&=f_1(w)+f_2(z)\\
&\le F(w)+F(z).
\end{aligned}
\]
Here the first inequality uses the supermodularity of $f_1$
and $\eta(x)\le\eta(z)$.
Thus, for fixed $a$, the function $F(a,\cdot)$ is supermodular.
Coordinatewise convexity follows because $F$ is the maximum
of two functions convex in each coordinate.

For fixed $A$, put $M=\min\{B(A),\vartheta\}$. For  $0\le a<\delta$, $F(a,A)$ varies with $a$ as $\kappa(a)M$, which has nonnegative first and second derivatives. On
$a>\delta$, its slope is $c$. Since the left-hand derivative at $a=\delta$ is
$c^2M\le c$, $F$ is convex and increasing in $a$. Moreover, for $a_1<a_2$, the increment $F(a_2,A)-F(a_1,A)$ depends on $A$ only through a nonnegative multiple of $M$, which is increasing in $A$. This establishes joint increasing directional convexity.

Finally, when $\alpha_p\le a$, the conditional density multiplier is at most $\kappa(a)$, so Lemma~\ref{lem:vector-comparison} bounds the conditional contribution of subsequent stages by $V_{\kappa(a)}(A)$, and \eqref{eq:joint-bound-branches} yields the stated upper bound.
\end{proof}

To apply Lemma~\ref{lem:joint-comparison} recursively over a finite family of admissible weights, we set $\Phi(\cdot)=0$ after the last prime coordinate and proceed by backward induction. At each step, if $z$ collects all remaining divisor weights, then $(a,A)=Lz$ for a matrix $L$ with nonnegative entries, and composition with $L$ preserves increasing directional convexity because
\[
\Psi(z+u+v)-\Psi(z+u)-\Psi(z+v)+\Psi(z)
=F(Lz+Lu+Lv)-F(Lz+Lu)-F(Lz+Lv)+F(Lz)\ge0
\]
for $\Psi(z)=F(Lz)$ and $z,u,v\ge0$. The function $\Delta(a)$ alone is concave on $(0,\delta)$, since $\Delta''(a)=-2(1-a)^{-3}<0$. Thus the full functional~\eqref{eq:joint-bound} is needed to maintain directional convexity.

Because both $V_{\kappa(a)}$ and $V_c$ are expectations of $\Phi$, the recursion is monotone in $\Phi$ and in the adjoined nonnegative weights. Choosing nested finite weight families with each $T_p$ replaced by $T_p\wedge b$ on a common probability space produces an increasing sequence of finite comparison bounds
\[0\le U_1\le U_2\le\cdots\le U_{\mathrm{sc}},\]
where $U_{\mathrm{sc}}$ is the sum of the reference stage bounds with the admissibility restrictions retained. Monotone convergence then yields the full comparison value $U=\lim_{n\to\infty}U_n\le U_{\mathrm{sc}}$.

\subsection{Accumulating decreases through a common positive kernel}
\label{sec:accumulated-reductions}

To obtain computable lower bounds for the reductions in the comparison
bound arising from~\eqref{eq:joint-bound}, let $Z_p$ denote the complete
admissible fibre weight at the prime $p$ and put
$\Delta_p(a)=c_p-\kappa_p(a)$. For the weight columns $A$, write
\[
 K_\rho f=\mathbb E_\rho f(S_T),\qquad
 \mathbb P_\rho(T\ge e)=\rho p^{-e}\quad(e\ge1),
\]
where the same random exponent $T$ is used in all columns. For finitely
many columns we use $S_{T\wedge b}$. Define the reference and reduced
comparison operators
\[
 \mathcal{C}_p f=g_p(Z_p)+K_{c_p}f,
\]
\begin{equation}\label{eq:reduced-comparison-operator}
 \mathcal{T}_p f = g_p(Z_p) + \max\{K_{\kappa_p(Z_p)} f, \, K_{c_p} f - \Delta_p(Z_p) \vartheta_p\}.
\end{equation}
Here each $\vartheta_p$ is fixed with $0\le \vartheta_p\le1/c_p$.

For unbounded $T_p$, define the common positive \emph{subprobability kernel}
\begin{equation}\label{eq:common-kernel}
M_p(0)=1-c_p/p,\qquad
M_p(v)=(p-1)p^{-v-1}\quad(v\ge1).
\end{equation}
Under either expectation kernel in~\eqref{eq:reduced-comparison-operator}, the probability of each value $v$ is at least $M_p(v)$. For the truncated exponent $T_p\wedge b$, the corresponding common kernel has the same coefficients for $0\le v<b$, coefficient $p^{-b}$ at $v=b$, and zero coefficients for $v>b$. In either case, a stronger state-dependent common subprobability kernel is obtained by keeping the coefficient at $v=0$ equal to $1-c_p/p$ and multiplying every coefficient at a positive value of $v$ by $\kappa_p(Z_p)$.

We also use $M_p$ to denote the associated operator:
\[
 M_p f=\sum_{v\ge0}M_p(v)f(S_v),
\]
with the corresponding finite sum over $0\le v\le b$ for the truncated kernel.

The following lemma uses a common kernel to give a lower bound for the difference between $\mathcal T_p f$ and $\mathcal T_p h$.
\begin{lemma}\label{lem:common-kernel-propagation}
If $f\ge h$ pointwise, then
\[
 \mathcal{T}_p f-\mathcal{T}_p h\ge M_p(f-h).
\]
The same assertion holds for the state-dependent stronger common kernel.
\end{lemma}

\begin{proof}
Select a branch in~\eqref{eq:reduced-comparison-operator} that attains the maximum for $\mathcal{T}_p h$, and let $K_\rho$ be its expectation kernel. Evaluating $\mathcal{T}_p f$ on the same branch gives
\[
\mathcal{T}_p f-\mathcal{T}_p h
\ge K_\rho(f-h)
\ge M_p(f-h).
\]
The second inequality follows from $f-h\ge0$ and the preceding lower bounds on the probabilities. No convexity assumption on $f-h$ is needed.
\end{proof}

We replace $\mathcal C_p$ by $\mathcal T_p$ in increasing prime order.
When the replacement at $p$ is made, all subsequent stages still use
the reference operators.
Let $B_p^0$ denote the quantity $B(A)$ defined in \eqref{eq:weight-column-differences}, with $\Phi$
taken to be the reference bound for the total contribution of
these subsequent stages.
By \eqref{eq:reduced-comparison-operator}, the local decrease is exactly
\begin{equation*}
    \Delta_p(Z_p)\min\{B_p^0,\vartheta_p\}.
\end{equation*}
We apply Lemma~\ref{lem:common-kernel-propagation} successively to the preceding comparison operators. Summing the resulting differences gives
\begin{equation}\label{eq:accumulated-reductions}
 U\le U_{\mathrm{sc}}
 -\sum_p\int \Delta_p(Z_p)\min\{B_p^0,\vartheta_p\}\,d\nu_{<p},
\end{equation}
where $\nu_{<p}$ is the subprobability measure obtained by applying the common subprobability kernels at the comparison primes preceding $p$, starting from the initial distribution of the comparison process.
When estimating the reductions over disjoint ranges of comparison primes, we may use a different valid common lower kernel for each range.

\section{Proof of Theorem~\ref{thm:odd-lcm}}\label{sec:odd-lcm}

We prove Theorem~\ref{thm:odd-lcm} by establishing $9\mid Q$
and $5\mid Q$ separately. We begin with a reduction that allows us
to assume that the modulus $3$ is present.
If $3\notin\mathcal D$, adjoining $0\pmod3$ gives an odd distinct
covering with least common multiple $Q'=\operatorname{lcm}(Q,3)$, and
\[
45\mid Q'\Longleftrightarrow45\mid Q.
\]
We may therefore assume throughout this section that $3\in\mathcal D$
and write $Q$ for the resulting least common multiple.

\subsection{Reduction modulo three and averaging}\label{sec:nine-reduction}

Suppose for contradiction that $9\nmid Q$. Then $v_3(Q)=1$.
Translate all residue classes so that the class with modulus $3$
is $0\pmod3$.
Any other class contained in it may be deleted.
The remaining integers lie in the two residue classes $1\pmod3$ and $2\pmod3$.

Because $v_3(Q)=1$, the integer $M=Q/3$ is coprime to $3$. Choose $s$ with $3s\equiv1\pmod M$. On the residue class indexed by $r\in\{1,2\}$ write $x=3t+r$ and introduce the coordinate
\begin{equation}\label{eq:affine-coordinate}
y=t+rs.
\end{equation}
The map in \eqref{eq:affine-coordinate} is a bijection on the product of all remaining cyclic coordinates.

If an original modulus is $m$ with $3\nmid m$, then
\begin{equation*}
3y=3t+3rs\equiv x\pmod m.
\end{equation*}
Thus the original class $a_m\pmod m$ becomes, on both residue classes, the same reduced class
\begin{equation*}
y\equiv sa_m\pmod m.
\end{equation*}
If the original modulus is $3m$ and its residue modulo $3$ is $r$, the same calculation gives
\begin{equation*}
y\equiv sa_{3m}\pmod m,
\end{equation*}
but this reduced class occurs only on the branch $x\equiv r\pmod 3$ that it intersects.

Let $\mathcal A$ be the union of the reduced classes from the original moduli coprime to $3$, and let $\mathcal B_r$ be the union of the reduced classes from the original moduli $3m$ whose residue modulo $3$ is $r$. Write $\mathcal D_{\mathcal{A}}$ and $\mathcal D_r$ for their respective sets of reduced moduli. The original covering gives
\begin{equation}\label{eq:reduced-covers}
\mathbb Z=\mathcal A\cup\mathcal B_1,
\qquad
\mathbb Z=\mathcal A\cup\mathcal B_2.
\end{equation}
Every reduced modulus exceeds $1$, since the original class of modulus
$3$ has been removed. Distinctness of the original moduli gives
\[
\mathcal D_1\cap\mathcal D_2=\varnothing.
\]
The reduced moduli are distinct within each of the three collections of classes. However,
$\mathcal D_{\mathcal{A}}\cap\mathcal D_r$ may be nonempty, since the same reduced modulus
$m$ can come from original moduli $m$ and $3m$.

We next compare the stage contributions of the two reduced coverings. Apply the distortion construction separately to the two coverings with
the same parameters $\delta_p$ and $c_p=(1-\delta_p)^{-1}$.
For $r=1,2$, let $B_{r,p}$ be the union of reduced classes introduced at
$p$, and let $\alpha_{r,p}$ be its proportion in each uniform new fibre.
With $\mu_{r,p}$ denoting the distorted measure, \eqref{eq:stage-measure}
gives
\[
L_{r,p}:=\mu_{r,p}(B_{r,p})
=c_p\mathbb E_{\mu_{r,p^-}}(\alpha_{r,p}-\delta_p)_+.
\]
A reduced modulus $m$ may arise from both original moduli $m$ and $3m$, but the distortion construction applies to the union of their classes.

Fix $p\ge5$.
Let $\mathcal S_{<p}$ be the set of positive integers supported on the primes preceding $p$.
For the common independent comparison integer $R$ on
these coordinates, define the divisor sums
\begin{equation*}
\begin{aligned}
S_{\mathcal{A},p}(R)&=\sum_{\substack{m\in\mathcal S_{<p},\ e\ge1\\mp^e\in\mathcal{D}_\mathcal{A}}}
 p^{-e}\mathbf{1}_{\{m\mid R\}},\\
X_{r,p}(R)&=\sum_{\substack{m\in\mathcal S_{<p},\ e\ge1\\mp^e\in\mathcal{D}_r}}
 p^{-e}\mathbf{1}_{\{m\mid R\}}\qquad(r=1,2).
\end{aligned}
\end{equation*}
Applying Theorem~\ref{thm:scalar-comparison} separately to each induced family gives
\begin{equation}\label{eq:labeled-branch-comparison}
L_{r,p}\le c_p\mathbb{E}\bigl(S_{\mathcal{A},p}(R)+X_{r,p}(R)-\delta_p\bigr)_+.
\end{equation}
Put
\[
H_p(R)=\sum_{m\mid R}\sum_{e\ge1}p^{-e}
=\frac{\tau(R)}{p-1}.
\]
Each pair $(m,e)$ occurs at most once in $S_{\mathcal{A},p}$, and
$\mathcal D_1\cap\mathcal D_2=\varnothing$ gives the separate bound
\[
S_{\mathcal{A},p}(R)\le H_p(R),
\qquad
X_{1,p}(R)+X_{2,p}(R)\le H_p(R).
\]
Lemma~\ref{lem:averaging}, with equal weights and
$\Phi(u)=(u-\delta_p)_+$, therefore gives, for each $R$,
\[
\frac12\sum_{r=1}^2
\bigl(S_{\mathcal{A},p}(R)+X_{r,p}(R)-\delta_p\bigr)_+
\le\frac12\bigl\{(H_p(R)-\delta_p)_+
                   +(2H_p(R)-\delta_p)_+\bigr\}.
\]
Taking expectations in the preceding inequality and applying \eqref{eq:labeled-branch-comparison}, with $X=\tau(R)$ and $a_p=\delta_p(p-1)$, we obtain
\begin{equation}\label{eq:averaged-stage-contribution}
\overline L_p:=\frac{L_{1,p}+L_{2,p}}2
\le\frac{c_p}{2(p-1)}
\left\{\mathbb E(X-a_p)_+
       +2\mathbb E(X-a_p/2)_+\right\}.
\end{equation}
Since both reduced systems are coverings, Corollary~\ref{cor:criterion} gives
\[
\sum_{p\ge5}\overline L_p
=\frac12\sum_{r=1}^2\sum_{p\ge5}L_{r,p}\ge1.
\]

\subsection{Finite estimates for the averaged contribution}\label{sec:nine-finite}

For $p\le7993$, take the integer $A_p$ from Table~\ref{tab:odd-distortion-parameters} in Appendix~\ref{app:parameters} and set $\delta_p=A_p/(p-1)$. Thereafter set
\begin{equation*}
\delta_p=\frac{71}{200}\quad(7993<p\le10^5),
\qquad
\delta_p=\frac12\quad(p>10^5).
\end{equation*}
All parameters lie in $[0,1/2]$.

During the finite stage $p\le10^5$, we retain upper bounds for $\mathbb{P}(\widetilde X=t)$
for $\widetilde X\le36000$, replacing any individual multiplier $D=T_p+1$
exceeding $d_*=10$ by $1$ as in \eqref{eq:replacement-bound}, so that
\begin{equation*}
\mathbb{P}(\widetilde D=1)=1-\frac{c_p}{p}+\frac{c_p}{p^{10}},
\qquad
\mathbb{P}(\widetilde D=d)=\frac{c_p(p-1)}{p^d}\quad(2\le d\le10).
\end{equation*}
Write $a_p=A/B$ with integers $A\ge0$ and $B>0$. If $\mathbb EX\le M^+/W$ and $\mathbb{P}(\widetilde X=t)\le P_t/W$, then \eqref{eq:replacement-bound} gives
\begin{equation}\label{eq:averaged-integer-bound}
\begin{aligned}
&\mathbb{E}(X-a_p)_+
 +2\mathbb{E}(X-a_p/2)_+\\
&\quad\le \frac{1}{BW}\biggl[
3BM^+-2AW
+\sum_{Bt<A}(A-Bt)P_t
+\sum_{2Bt<A}(A-2Bt)P_t
\biggr].
\end{aligned}
\end{equation}
Since $a_p<35500$ for $p\le10^5$, the cutoff $H=36000$ covers both sums in~\eqref{eq:averaged-integer-bound}.

For $z\ge0$, \eqref{eq:cubic} gives
\begin{equation}\label{eq:averaged-cubic-bound}
\begin{aligned}
(z-1/2)_++(2z-1/2)_+
&=\max\{2(z-1/4)_+,\,3(z-1/3)_+\}\\
&\le\max\{(128/27)z^3,\,4z^3\}
=\frac{128}{27}z^3.
\end{aligned}
\end{equation}
For $p>10^5$, substitute $z=X/(p-1)$ and $c_p=2$ into
\eqref{eq:averaged-stage-contribution} to obtain
\[
\overline L_p\le\frac{128}{27}
\frac{\mathbb E X^3}{(p-1)^3}.
\]

\begin{proposition}\label{prop:nine-numerics}
With the parameters and recursions above, integer arithmetic with directed rounding at scale $10^{18}$ gives
\begin{equation}\label{eq:nine-bounds}
\begin{aligned}
\sum_{5\le p\le10^5}\overline L_p
&\le\NineFiniteUpper,\\
\sum_{10^5<p\le10^7}\overline L_p
&\le\NineMiddleUpper,\\
\sum_{p>10^7}\overline L_p
&\le\NineTailUpper.
\end{aligned}
\end{equation}
\end{proposition}

\begin{proof}
Put $W=10^{18}$, and write $M_1=\mathbb E X$ and
$M_3=\mathbb E X^3$ for the full moments before the current prime.
The initial values are
\[
P_1=W,\qquad P_t=0\ (t>1),\qquad M_1=M_3=1.
\]
For $p\le10^5$, compute the stage bound from~\eqref{eq:averaged-integer-bound}
before updating the probability numerators by~\eqref{eq:product-recursion}.
Only values $t=2^a3^b5^c7^d\le36000$, with $a,b,c,d\ge0$, are needed,
since every retained multiplier is at most $10$.
The untruncated moment updates are
\[
M_1'=M_1\left(1+\frac{c_p}{p-1}\right),\qquad
M_3'=M_3\left(1+c_p\frac{7p^2-2p+1}{(p-1)^3}\right).
\]
All scaled upper bounds are rounded upward. For $10^5<p\le10^7$,
use~\eqref{eq:averaged-cubic-bound} for the stage bound and continue the full
third-moment update. Since
\[
1+\frac{2(7p^2-2p+1)}{(p-1)^3}
\le1+\frac{15}{p-1},
\]
Lemma~\ref{lem:log-tail} applies with $\beta=15$.
Put $J=664579$, the number of primes up to $10^7$, and
$b=J-1=664578$.
The series for $e$ and an integer comparison give
\[
e<2+\frac{1/2}{1-1/3}=\frac{11}{4},\qquad
11^{13}<664578\cdot4^{13},
\]
so $e^{13}<(11/4)^{13}<b$ and $\log b>13>(15-3)/2$.
If $\mathcal M_3$ bounds the full third moment after $p_J$, the lemma
with $\ell=13$ gives
\[
\sum_{p>10^7}\overline L_p
\le\frac{128\mathcal M_3}{27\cdot13^2\cdot14b^2}.
\]
Evaluating these upward-rounded bounds over the three ranges yields
\eqref{eq:nine-bounds}.
\end{proof}

Adding the three bounds in Proposition~\ref{prop:nine-numerics} gives
\begin{equation*}
\sum_{p\ge5}\overline L_p
\le\NineTotalUpper<1,
\end{equation*}
contradicting the existence of both coverings in \eqref{eq:reduced-covers}. Thus $9\mid Q$.

\subsection{Divisibility by five}\label{sec:five-divisibility}

Suppose for contradiction that $5\nmid Q$. Let $S$ be the
complement of the union of the classes whose moduli are powers of
$3$, and condition the uniform measure on $S$, as in
Lemma~\ref{lem:pure-three-conditioning}.
By Corollary~\ref{cor:conditioned-distortion}, the $3$-coordinate is
compared with a nonnegative integer-valued random exponent $T_3$ satisfying
\[
\mathbb P(T_3\ge e)=\frac2{3^e}\qquad(e\ge1).
\]
Its divisor-count multiplier $D_3=T_3+1$ satisfies
\begin{equation*}
\begin{gathered}
\mathbb P(D_3=1)=\frac13,\qquad
\mathbb P(D_3=d)=\frac4{3^d}\quad(d\ge2),\\
\mathbb E D_3=2,\qquad \mathbb E D_3^3=\frac{31}{2}.
\end{gathered}
\end{equation*}
These moments match~\eqref{eq:moments} with $q=3$ and $\gamma_3=2$.

We process the remaining prime coordinates $p\ge7$ in increasing order with $\delta_p=3/10$ and
$c_p=10/7$ throughout. Let $X$ be the product of the divisor-count
multipliers preceding $p$, including $D_3$.
By \eqref{eq:divisor-upper} and \eqref{eq:stage-measure},
\begin{equation}\label{eq:no-five-contribution}
L_p\le\frac{10}{7(p-1)}
\mathbb{E}\left(X-\frac{3(p-1)}{10}\right)_+.
\end{equation}
For $p\le5000$, we retain the distribution of $\widetilde X$ up to $H=1600$ with multiplier cutoff $d_*=10$ as in~\eqref{eq:replacement-bound}, while updating the untruncated first and third moments via~\eqref{eq:moments}. For $5000<p\le10^7$, applying~\eqref{eq:cubic} to~\eqref{eq:no-five-contribution} bounds $L_p$ by $(4000/1701)\mathbb{E}X^3/(p-1)^3$, and Lemma~\ref{lem:log-tail} with $\beta=15$ and $\ell=13$ bounds the tail $p>10^7$. Integer arithmetic with directed rounding at scale $10^{18}$ yields
\begin{equation*}
\begin{aligned}
\sum_{7\le p\le5000}L_p&\le\FiveFiniteUpper,\\
\sum_{5000<p\le10^7}L_p&\le\FiveMiddleUpper,\\
\sum_{p>10^7}L_p&\le\FiveTailUpper,
\end{aligned}
\end{equation*}
whose sum is
\begin{equation*}
\sum_{p\ge7}L_p\le\FiveTotalUpper<1.
\end{equation*}
By Corollary~\ref{cor:conditioned-distortion}, the remaining classes cannot cover $S$, so $5\mid Q$. Together with $9\mid Q$, this proves $45\mid Q$ and completes the proof of Theorem~\ref{thm:odd-lcm}.

\section{Proof of Theorem~\ref{thm:odd}}\label{sec:odd}
We prove Theorem~\ref{thm:odd} in two steps. We first rule out
the possibility that every modulus is at least $15$. We then
exclude the remaining cases in which the least modulus is $11$
or $13$.
\subsection{Ruling out minimum modulus at least fifteen
}\label{sec:odd-fifteen}

We first describe the finite recursion used to bound the stage contributions under the assumption that every modulus is at least $15$.
Suppose that every modulus is an odd integer at least $15$. For the current prime $p$, let $e_0$ be the least positive integer such that $p^{e_0}\ge15$. For the comparison integer $R$ supported on the preceding primes, define
\begin{equation}\label{eq:odd-threshold}
Z_p(R)=
\sum_{1\le e<e_0}
\frac{\tau(R)-D_{\left\lceil 15/p^e\right\rceil}(R)}{p^e}
+\frac{\tau(R)}{p^{e_0-1}(p-1)},
\end{equation}
where the final term sums all exponents $e\ge e_0$.

Only divisors of $R$ below $15$ enter the subtracted terms in \eqref{eq:odd-threshold}. They are determined by the exponent tuple
\begin{equation*}
\sigma=\bigl(\min\{T_3,2\},\min\{T_5,1\},\min\{T_7,1\},
\min\{T_{11},1\},\min\{T_{13},1\}\bigr),
\end{equation*}
which has $3\cdot2^4=48$ possible values. Before a prime $q$ is processed, set $\sigma_q=0$. Define $R(\sigma)=3^{\sigma_3}5^{\sigma_5}7^{\sigma_7}11^{\sigma_{11}}13^{\sigma_{13}}$. Then $R(\sigma)\mid R$, and hence
\[D_a(R(\sigma))\le D_a(R),\]
so replacing the latter by the former in the subtraction terms of \eqref{eq:odd-threshold} remains a valid upper bound.

The prime order need not be increasing, as in
Section~\ref{sec:distortion}. We use the five initial steps:
\[
\begin{array}{c|ccccc}
p&3&5&13&7&11\\ \hline
\delta_p&11/200&161/3200&311/3200&303/1600&333/1600
\end{array}
\]
All subsequent primes are processed in increasing order with the
parameters in Table~\ref{tab:odd-distortion-parameters} and
Appendix~\ref{app:parameters}.

Put $H=36000$ and $W=10^{18}$, and let $\sigma(R)$ denote the tuple
above. Before each initial prime step, we retain scaled integer upper bounds
\[
P(t,\sigma)\ge W\cdot\mathbb P\bigl(\tau(R)=t,\ \sigma(R)=\sigma\bigr)
\qquad(1\le t\le H),
\]
starting from $P(1,0)=W$ and zero elsewhere.
Here $Z_p(t,\sigma)$ denotes
\eqref{eq:odd-threshold} with $\tau(R)=t$ and each $D_a(R)$ replaced
by $D_a(R(\sigma))$. Since $Z_p(t,\sigma)\le t/(p-1)$,
\begin{equation}\label{eq:odd-truncation}
\mathbb E\left[\frac{c_p\tau(R)}{p-1}
                \mathbf1_{\{\tau(R)>H\}}\right]
\le\frac{c_p\mathbb E\tau(R)^3}{(p-1)H^2}.
\end{equation}
Thus the initial stage contribution is bounded by
\[
L_p=\mu_p(B_p)\le
\frac{c_p}{W}\sum_{\substack{t\le H\\\sigma}}
P(t,\sigma)(Z_p(t,\sigma)-\delta_p)_+
+\frac{c_p\mathbb E\tau(R)^3}{(p-1)H^2}.
\]
We evaluate this bound before each step. For $T_p=v$, the pair $(t,\sigma)$ is updated according to
\[
(t,\sigma)\longmapsto(t(v+1),\sigma'),\qquad
\sigma'_p=\begin{cases}\min\{v,2\},&p=3,\\
                         \min\{v,1\},&p\ne3,
             \end{cases}
\]
with $\sigma'_r=\sigma_r$ for $r\ne p$.
For integer bounds $\Pi_p(v)\ge W\cdot\mathbb P(T_p=v)$, we set
\[
P'(u,\eta)=\sum_{\substack{t(v+1)=u\\\sigma'=\eta}}
\left\lceil\frac{P(t,\sigma)\Pi_p(v)}W\right\rceil
\qquad(1\le u\le H),
\]
where the sum runs over preceding pairs $(t,\sigma)$ and $v\ge0$. Every multiplier is at least $1$, so pairs with $t>H$ cannot contribute to $u\le H$ during these five steps. We update the full first and third moments independently by~\eqref{eq:moments}.

After these five primes, every remaining prime satisfies $p\ge17$, so
\eqref{eq:odd-threshold} reduces to
\[Z_p(R)=\tau(R)/(p-1).\]
Let $X_0$ be the full, untruncated divisor count after the five steps,
and set $Y_0=X_0$. If $P_0(t, \sigma)$ denotes the resulting scaled upper bounds after these five steps, put
\[
Q_0(t)=\sum_\sigma P_0(t,\sigma)
\ge W\cdot\mathbb P(X_0=t)\qquad(1\le t\le H).
\]
For the subsequent primes $p\le10^5$, apply
Lemma~\ref{lem:whole-reset} with $D_j=T_p+1$ and $d_*=10$. The probabilities are
\[
\pi_j(1)=1-\frac{c_p}{p},\qquad
\pi_j(d)=\frac{c_p(p-1)}{p^d}\quad(2\le d\le10),\qquad
\rho_j=\frac{c_p}{p^{10}}.
\]
The tail term at $u = 1$ accounts for all values of $Y_{j-1}$, since
\[
\mathbb P(Y_j=1,\ D_j>10)
=\rho_j\sum_{u\ge1}\mathbb P(Y_{j-1}=u)=\rho_j,
\]
including those with $Y_{j-1}>H$.
Put $a=\delta_p(p-1)<H$.
By~\eqref{eq:whole-reset-positive} and
$q_{j-1}(u)\le Q_{j-1}(u)/W$,
\[
\begin{aligned}
L_p=\mu_p(B_p)
&\le\frac{c_p}{p-1}\mathbb E(X_{j-1}-a)_+\\
&\le\frac{c_p}{p-1}\left\{
\mathbb E X_{j-1}-a+
\frac1W\sum_{1\le u<a}(a-u)Q_{j-1}(u)\right\}.
\end{aligned}
\]
After evaluating the stage bound, we update $Q_j$ by \eqref{eq:whole-reset-integer} and the full moments by~\eqref{eq:moments}.

\begin{proposition}\label{prop:odd-fifteen}
No finite family of residue classes with distinct odd moduli at least
$15$ covers $\mathbb Z$.
\end{proposition}

\begin{proof}
Suppose that such a family covers $\mathbb Z$. The initial computation, including the errors in
\eqref{eq:odd-truncation}, gives
\[
\sum_{p\in\{3,5,13,7,11\}}L_p\le\OddFifteenEarlyUpper
\]
and
\[
\mathbb E X_0\le\OddFifteenFirstMomentUpper,\qquad
\mathbb E X_0^3\le\OddFifteenThirdMomentUpper.
\]
We continue the finite recursion through $10^5$ and set $\delta_p=1/2$ for $p>10^5$. This choice gives the cubic coefficient $32/27$. On $10^5<p\le10^7$, we update the third moment using $1+14.001/(p-1)$ as an upper bound for the third-moment multiplier, which is valid since
\begin{equation}\label{eq:odd-moment-inflation}
\frac{2(7p^2-2p+1)}{(p-1)^2}
=14+\frac{24p-12}{(p-1)^2}<14.001\qquad(p>10^5).
\end{equation}
With $J=664579$ and $\mathcal M_3$ denoting the third-moment upper bound after $p_J\le10^7$, we apply Lemma~\ref{lem:log-tail} with $\beta=15$, $\ell=13$, $b=J-1$, and $C=32/27$ to obtain
\begin{equation}\label{eq:odd-tail-relaxation}
\sum_{p>10^7}L_p
\le\frac{32\mathcal M_3}{27\cdot2366(J-1)^2}
\le\frac{32\mathcal M_3}{27\cdot1870J^2}.
\end{equation}
The bounds for the four ranges are summarized below.
\begin{center}
\begin{tabular}{lr}
\toprule
Range & Upper bound for the sum of contributions\\\midrule
Five initial primes & $\OddFifteenEarlyUpper$\\
Remaining $p\le10^5$ & $\OddFifteenFiniteUpper$\\
$10^5<p\le10^7$ & $\OddFifteenMiddleUpper$\\
$p>10^7$ & $\OddFifteenTailUpper$\\\midrule
Total & $\OddFifteenTotalUpper$\\
\bottomrule
\end{tabular}
\end{center}
The total is less than $1$, so Corollary~\ref{cor:criterion}
shows that these residue classes cannot form a covering.
\end{proof}

\subsection{The cases of least modulus eleven and thirteen}\label{sec:odd-eleven-thirteen}

By Proposition~\ref{prop:odd-fifteen}, a least odd modulus exceeding
$9$ must be $q\in\{11,13\}$. Translate its class to $0\pmod q$ and
equip the $q$-coordinate with the uniform reference measure $\nu_q$ on
its complement, leaving the other coordinates uniform.
Every residue class modulo $q^e$ then has $\nu_q$-measure at most
$q^{1-e}/(q-1)$. Lemma~\ref{lem:reference-law} gives the coefficient
\begin{equation}\label{eq:conditioned-q}
\widehat c_q=\frac q{q-1}\frac1{1-\delta_q},\qquad
\mathbb P(T_q\ge e)=\widehat c_q q^{-e}\quad(e\ge1),
\end{equation}
for the comparison exponent $T_q$, while the stage factor remains
$c_q=(1-\delta_q)^{-1}$.

Write $Z_{p,q}(R)$ for the comparison fibre weight at $p$ in this case.
At $p=q$, removing the excluded modulus $q$, corresponding to $m=1$
and $e=1$, gives
\begin{equation*}
Z_{q,q}(R)=\frac q{q-1}
\left(\frac{\tau(R)}{q-1}-\frac1q\right)
=\frac{q\tau(R)-(q-1)}{(q-1)^2}.
\end{equation*}
At $p\ne q$, a pair $(m,e)$ is excluded precisely when the full modulus
$mp^e<q$, equivalently $m<\lceil q/p^e\rceil$. Hence
\[
Z_{p,q}(R)=\sum_{e\ge1}p^{-e}
\left\{\tau(R)-D_{\lceil q/p^e\rceil}(R)\right\}
\qquad(p\ne q).
\]
For example, when $q=11$ and $p=7$, the term $m=3$, $e=1$ is retained
whenever $3\mid R$, since its full modulus is $21$.

For $q=11$, the four initial steps are
\[
\begin{array}{c|cccc}
p&3&5&11&7\\ \hline
\delta_p&11/200&1/20&77/640&19/100
\end{array}
\]
For $q=13$, use the five initial primes $3,5,13,7,11$ and their
parameters from Subsection~\ref{sec:odd-fifteen}.
In both cases, process the remaining primes in increasing order using
Table~\ref{tab:odd-distortion-parameters}. These primes satisfy $p>q$, so
$Z_{p,q}(R)=\tau(R)/(p-1)$.

\begin{proposition}\label{prop:odd-eleven-thirteen}
The least modulus of an odd distinct covering system cannot be $11$ or $13$.
\end{proposition}

\begin{proof}
We use the comparison exponent coefficient~\eqref{eq:conditioned-q} at $q$ and the weights $Z_{p,q}$ above in the initial recursion for $(t,\sigma)$. The terms with $\tau(R)>H=36000$ contribute at most
\[
s_{p,q}\cdot\frac{\mathbb E\tau(R)^3}{H^2},\qquad
s_{p,q}=\begin{cases}
\displaystyle\frac{q c_q}{(q-1)^2},&p=q,\\[5pt]
\displaystyle\frac{c_p}{p-1},&p\ne q.
\end{cases}
\]
This follows from $c_p(Z_{p,q}-\delta_p)_+\le s_{p,q}\tau(R)$ and
$\tau(R)\mathbf1_{\{\tau(R)>H\}}\le\tau(R)^3/H^2$.
After the four initial steps for $q=11$, or five for $q=13$, we sum the upper bounds $P_0(t, \sigma)$ over $\sigma$ and use the full untruncated moments as in Subsection~\ref{sec:odd-fifteen}. We apply Lemma~\ref{lem:whole-reset} with $H=36000$ and $d_*=10$,
including the tail term $\mathbf{1}_{\{u=1\}}\rho_j$. The resulting finite stage bounds, together with \eqref{eq:odd-moment-inflation} and~\eqref{eq:odd-tail-relaxation}, give
\begin{center}
\begin{tabular}{lrr}
\toprule
Range & $q=11$ & $q=13$\\\midrule
Initial primes & $\OddElevenEarlyUpper$ & $\OddThirteenEarlyUpper$\\
Remaining $p\le10^5$ & $\OddElevenFiniteUpper$ & $\OddThirteenFiniteUpper$\\
$10^5<p\le10^7$ & $\OddElevenMiddleUpper$ & $\OddThirteenMiddleUpper$\\
$p>10^7$ & $\OddElevenTailUpper$ & $\OddThirteenTailUpper$\\\midrule
Total & $\OddElevenTotalUpper$ & $\OddThirteenTotalUpper$\\
\bottomrule
\end{tabular}
\end{center}
Since both totals are strictly less than $1$, Corollary~\ref{cor:criterion} rules out a covering of the complement of $0\pmod q$ in either case.
\end{proof}

Propositions \ref{prop:odd-fifteen} and \ref{prop:odd-eleven-thirteen} together show that the least odd modulus is at most $9$, completing the proof of Theorem~\ref{thm:odd}.

\section{Proof of Theorem~\ref{thm:bounds}}\label{sec:bounds-proof}
We first prove the bound for general moduli and then turn to the squarefree case.
\subsection{General moduli}\label{sec:general}

Suppose for contradiction that every modulus satisfies $d\ge K:=14000$.

\medskip\noindent\emph{Admissible fibre weights and the initial union.}
We equip the coordinates corresponding to
\[
\mathcal P_0=\{2,3,5,7,11,13,17,19,23\}
\]
with the product uniform measure and begin the distortion process at $p=29$. For $p\ge29$, put
$R_{<p}=\prod_{r<p}r^{T_r}$, and use the local abbreviations
$h_{p,1}=\lceil K/p\rceil$ and $h_{p,2}=\lceil K/p^2\rceil$.
Since $p^3>K$, the complete admissible fibre weight is
\[
 Z_p=\frac{\tau(R_{<p})}{p-1}-C_p,\qquad
 C_p=\frac{D_{h_{p,1}}(R_{<p})}{p}
                 +\frac{D_{h_{p,2}}(R_{<p})}{p^2}\ge0.
\]
Here $D_a(n)$, as defined in \eqref{eq:small-divisors-and-support}, denotes the number of positive divisors of $n$ that are less than $a$.
Let
\[
 L_p^{\mathrm{sc}}=c_p\mathbb{E}\left(\frac{\tau(R_{<p})}{p-1}-\delta_p\right)_+
\]
and
\begin{equation}\label{eq:small-divisor-reduction}
 G_p=c_p\mathbb{E}\min\left\{C_p,
                  \left(\frac{\tau(R_{<p})}{p-1}-\delta_p\right)_+\right\}.
\end{equation}
Then $\mathbb{E}g_p(Z_p)=L_p^{\mathrm{sc}}-G_p$, with $G_p=0$ for $p\ge K$. Thus the sum of the reference stage bounds is
\[
U_{\mathrm{sc}}
=\sum_{p\ge29}\mathbb{E}g_p(Z_p)
=\sum_{p\ge29}\bigl(L_p^{\mathrm{sc}}-G_p\bigr).
\]

For each nonempty support $S\subseteq\mathcal P_0$, the probability of the union of initial classes with support $S$ is at most
\[
 w_S=\prod_{p\in S}\frac1{p-1}
 -\sum_{\substack{2\le d<K\\\operatorname{supp}(d)=S}}\frac1d.
\]
The events associated with disjoint supports are independent. For $U\subseteq\mathcal P_0$, compute
\[
 Q_\varnothing=1,\qquad
 Q_U=Q_{U\setminus\{b\}}
        -\sum_{\substack{S\subseteq U\\b\in S}}w_SQ_{U\setminus S}
       \quad(b\in U).
\]
All $512$ values $Q_U$ are positive exact fractions. To check the hypotheses of Lemma~\ref{lem:shearer} for every subfamily of the $511$ support events, let $Q_U(u)$ denote the same polynomial with weights $u_S\in[0,w_S]$. For nonempty $S\subseteq U$,
\[
 \frac{\partial Q_U(u)}{\partial u_S}=-Q_{U\setminus S}(u).
\]
Starting from $Q_\varnothing=1$, induction on $|U|$ gives
\[
Q_U(u)\ge Q_U(w)>0
\]
throughout the box. Setting the weights of omitted support events to zero verifies positivity for every subfamily. Lemma~\ref{lem:shearer} now bounds the initial union by
\begin{equation}\label{eq:initial-union-bound}
 I_0=1-Q_{\mathcal P_0}
 =\frac{183094911650617622458764685359451238903}
 {1764570451786702440402513892265164800000},
\end{equation}
and marginal preservation keeps this bound unchanged at subsequent stages.

\medskip\noindent\emph{Lower bound for the small-divisor reduction.}
The initial factor $R_{<29}$ has distribution
\begin{equation}\label{eq:initial-comparison-law}
 \mathbb{P}(R_{<29}=n)=\frac{110592}{676039}\cdot\frac1n
       \quad(n\text{ is }23\text{-smooth}).
\end{equation}
For $n\le10^{12}$, the probabilities rounded down at scale $10^{18}$ sum to $0.999999954168606443$. Omitting larger $n$ reduces $G_p$. The reference moments remain untruncated.

For $29\le q\le499$, we have $h_{q,1}\le483<29^2$. A divisor counted by $D_{h_{q,1}}(R_{<q})$ therefore contains at most one prime $r\ge29$, to the first power, and
\[
 C_q=\frac{D_{h_{q,1}}(n)+\Lambda}{q}+\frac{D_{h_{q,2}}(n)}{q^2},\qquad
 \Lambda=\sum_{29\le r<q}
 D_{\lceil K/(qr)\rceil}(n)\mathbf1_{\{T_r\ge1\}}.
\]
The coefficients of the indicators in $\Lambda$ involve only divisors at most $15$ and are determined by the $4\cdot3\cdot2^4=192$ tuples
\[
 \left(\min\{v_2(n),3\},\min\{v_3(n),2\},
 \mathbf1_{\{5\mid n\}},\mathbf1_{\{7\mid n\}},
 \mathbf1_{\{11\mid n\}},\mathbf1_{\{13\mid n\}}\right).
\]
Put $\overline D_r=\min\{D_r,10\}$, where the random multipliers $D_r=T_r+1$ are those of Lemma~\ref{lem:vector-comparison}. Its full distribution is
\[
 \mathbb P(\overline D_r=d)=
 \begin{cases}
 1-c_r/r,&d=1,\\
 c_r(r-1)/r^d,&2\le d\le9,\\
 c_r/r^9,&d=10,
 \end{cases}
 \qquad
 \mathbf1_{\{\overline D_r\ge2\}}=\mathbf1_{\{T_r\ge1\}}.
\]
We compute the joint distribution of
\[
 \left(\min\left\{512,\prod_{29\le r<q}\overline D_r\right\},
       \Lambda\right).
\]
Replacing $D_r$ by $\overline D_r$ can only lower the integrand in~\eqref{eq:small-divisor-reduction}. Replacing a product above $512$ by $512$ leaves it unchanged, since $\delta_q\le1/2$ and
\[
 C_q\le\frac{h_{q,1}-1}{q}+\frac{h_{q,2}-1}{q^2}
       <\frac{512}{q-1}-\frac12.
\]
Thus downward rounding gives a lower bound, with all larger product values represented by $512$.

For $499<q<K$, we have $h_{q,1}\le28$ and $h_{q,2}=1$, so $C_q=D_{h_{q,1}}(n)/q$. We use the one-dimensional distribution of
\[
 \min\left\{8192,\prod_{29\le r<q}\overline D_r\right\}.
\]
Here $C_q\le(h_{q,1}-1)/q<8192/(q-1)-1/2$, so the same argument applies. Downward rounding gives
\begin{equation}\label{eq:small-divisor-reduction-bound}
 \sum_{29\le p<K}G_p\ge
              \frac{1230344533175833263}{10^{18}}.
\end{equation}

\medskip\noindent\emph{Reductions from smaller conditional density multipliers.}
Fix $q>p\ge29$ and condition on all comparison coordinates in $[2,q)\setminus\{p\}$. Put
\[
 R^*=\frac{R_{<q}}{p^{T_p}},\qquad b=\frac{\tau(R^*)}{q-1},
\]
\[
 C_0=\frac{D_{h_{q,1}}(R^*)}{q}+\frac{D_{h_{q,2}}(R^*)}{q^2},
 \qquad C_1=\frac{D_{\lceil K/(pq)\rceil}(R^*)}{q}.
\]
Since $p^2q,pq^2,q^3>K$, we have $C_q=C_0+\mathbf1_{\{T_p\ge1\}}C_1$. Also $C_0\le(1-q^{-2})b<b$ and $C_1\le(1-q^{-1})b<b$.
The nonnegative weight columns are $b-C_0,b-C_1,b,b,\ldots$, with partial sums
\[
 S_v=(v+1)b-C_0-\mathbf1_{\{v\ge1\}}C_1\qquad(v\ge0).
\]
In particular, $Z_q=S_v$ when $T_p=v$. Its contribution to the weighted first difference at $p$ is
\begin{equation}\label{eq:local-weighted-difference}
 B_{pq}=\sum_{v\ge1}p^{-v}\{g_q(S_v)-g_q(S_{v-1})\}\ge0.
\end{equation}
Recall from Section~\ref{sec:accumulated-reductions} that $B_p^0$
is obtained from $B(A)$ by taking $\Phi$ to be the reference
bound for the total contribution of subsequent stages. For any finite set of subsequent primes,
\[
 B_p^0\ge\sum_q\mathbb E[B_{pq}\mid R_{<p}],
\]
where the conditional expectations are taken over the independent exponents $T_r$ with $p<r<q$, each satisfying $\mathbb P(T_r\ge e)=c_r r^{-e}$ for $e\ge1$, while $R_{<p}$ is held fixed.
We denote the reduction sums in~\eqref{eq:accumulated-reductions} over $29\le p\le233$ and $233<p\le63541$ by $\mathcal R_{[29,233]}$ and $\mathcal R_{>233}$, respectively.

For $29\le p\le233$, we retain $23$-smooth integers $n\le10^{12}$ from~\eqref{eq:initial-comparison-law} with $\tau(n)\le65$. Let $w$ be the exact product of the preceding factors $D_r$, $29\le r<p$, under the common kernel~\eqref{eq:common-kernel}, and retain values with $h=\tau(n)w\le65$. Omitting other values lowers the nonnegative reduction integral. Since $n\mid R_{<p}$ and $\Delta_p$ is nonincreasing, retaining only the exclusions from $\mathcal P_0$ gives
\[
 \Delta_p(Z_p)\ge
 \Delta_p\!\left(\frac h{p-1}
       -\frac{D_{h_{p,1}}(n)}p-\frac{D_{h_{p,2}}(n)}{p^2}\right).
\]
For $p<q\le491$, consider the list
\[
 \left(D_{\lceil K/(qr)\rceil}(n)\right)_{\substack{29\le r<p\\r\text{ prime}}},
\]
with one entry for each prime.
Let $W_{pq}(n,w)$ be the sum of its largest $\Omega(w)$ entries, or of all its entries if the list has fewer than $\Omega(w)$ entries.
At most $\Omega(w)$ of these primes have $T_r\ge1$, because each such prime contributes a factor $T_r+1\ge2$ to $w$. With
\[
 \Lambda_{pq}=\sum_{p<r<q}D_{\lceil K/(qr)\rceil}(n)\mathbf1_{\{T_r\ge1\}},
\]
we obtain
\[
 \overline C_0=\frac{D_{h_{q,1}}(n)+W_{pq}(n,w)+\Lambda_{pq}}q
              +\frac{D_{h_{q,2}}(n)}{q^2}\ge C_0,
 \qquad \underline b=\frac h{q-1}\min\left\{512,\prod_{p<r<q}\overline D_r\right\}\le b.
\]
Here $C_1=D_{\lceil K/(pq)\rceil}(n)/q\le\underline b$. The $v=1$ term of~\eqref{eq:local-weighted-difference} gives
\[
\begin{aligned}
 B_{pq}&\ge\frac{g_q(2b-C_0-C_1)-g_q(b-C_0)}p\\
 &=\frac{c_q}p\min\{b-C_1,(2b-C_0-C_1-\delta_q)_+\}\\
 &\ge\frac{c_q}p\min\{\underline b-C_1,
           (2\underline b-\overline C_0-C_1-\delta_q)_+\}.
\end{aligned}
\]
The last inequality holds because the preceding expression is nondecreasing in $b$ and nonincreasing in $C_0$.
We keep $\Lambda_{pq}$ exact because increasing it lowers this bound.
We add the expectations for all $p<q\le491$ to those from the one-dimensional recursion for $499\le q\le63541$, where $C_1=0$, before taking the minimum with $\vartheta_p=1/c_p$. This yields
\begin{equation}\label{eq:early-reduction-bound}
 \mathcal R_{[29,233]}\ge\frac{8837268985393504}{10^{18}}.
\end{equation}

Put
\[
 C_p^{\max}=\frac{h_{p,1}-1}{p}+\frac{h_{p,2}-1}{p^2}\qquad(p\ge29).
\]
For $233<p\le63541$ and $h=\tau(R_{<p})$, we have $p^2>K$, so
\[
 \max\{0,h/(p-1)-C_p^{\max}\}\le Z_p
       \le\frac h{p-1}-\frac{\mathbf1_{\{p<K\}}}{p}.
\]
We evaluate $\kappa_p$ at the lower bound to strengthen the common kernel, and $\Delta_p$ at the upper bound to lower the local reduction.

Since $pq>K$ for $q>p>233$, we have $C_1=0$. Each increment
$g_q((v+1)b-C_0)-g_q(vb-C_0)$ decreases in $C_0$. Replacing $C_0$ by $C_q^{\max}$ therefore gives a lower bound for $B_{pq}$. For fixed $p$, put
\[
 V_{pq}=\prod_{p<r<q}D_r,\qquad
 \theta_q=(q-1)(\delta_q+C_q^{\max}),
\]
\[
\begin{aligned}
 f(h)&=\sum_{p<q\le63541}\frac{c_q}{q-1}
                    \mathbb E(hV_{pq}-\theta_q)_+,\\
 \alpha_{\mathrm{lin}}&=\sum_{p<q\le63541}\frac{c_q}{q-1}\mathbb EV_{pq},\\
 u(h)&=\sum_{p<q\le63541}\frac{c_q}{q-1}
                         \mathbb E(\theta_q-hV_{pq})_+.
\end{aligned}
\]
These expectations are taken with respect to the reference probability measure. Since $V_{pq}\ge1$ and
$\max_{29\le q\le63541}\theta_q=1304720829/50000<26095$, we have $u(h)=0$ for $h\ge26095$. Applying \eqref{eq:positive-part} cancels the terms independent of $h$ in the differences, so
\begin{equation}\label{eq:scalar-weighted-difference}
\begin{aligned}
 \sum_{v\ge1}p^{-v}\{f((v+1)h)-f(vh)\}
 &=\frac{\alpha_{\mathrm{lin}}h}{p-1}
       +\sum_{v\ge1}p^{-v}\{u((v+1)h)-u(vh)\}\\
 &=\frac{\alpha_{\mathrm{lin}}h}{p-1}-\frac{u(h)}p
       +\sum_{d\ge2}\frac{p-1}{p^d}u(dh).
\end{aligned}
\end{equation}
Only finitely many terms involving $u$ are nonzero. We use lower bounds for $\alpha_{\mathrm{lin}}$ and the positive terms, and an upper bound for the subtracted term $u(h)$. In the backward recursion for $u$, an upper bound on an expectation of the next function is obtained from
\[
 \mathbb E u_{\mathrm{next}}(hD_r)
 \le\sum_{d=1}^{12}\mathbb P(D_r=d)u^+_{\mathrm{next}}(hd)
      +\frac{c_r}{r^{12}}\sup_{t\ge1}u^+_{\mathrm{next}}(t),
\]
where $u^+_{\mathrm{next}}$ is an upper bound obtained by upward rounding. The last term includes the entire probability $\mathbb P(D_r>12)=c_r/r^{12}$. Lower bounds omit this nonnegative tail. The terms with $d>12$ in the last sum of~\eqref{eq:scalar-weighted-difference} may likewise be omitted in a lower bound. Arithmetic with directed rounding gives
\begin{equation}\label{eq:late-reduction-bound}
 \mathcal R_{>233}\ge\frac{5740271089542}{2^{50}}.
\end{equation}

\medskip\noindent\emph{Scalar upper bound and completion of the proof.}
The $6360$ parameters $\delta_p$ for $29\le p\le63541$ are fixed rationals with denominator $10^6$, as specified in Appendix~\ref{app:parameters}. We then take
\[
\delta_p=\begin{cases}
41323/100000&\text{if } 63541<p\le10^7,\\
1/2&\text{if } p>10^7.
\end{cases}
\]
For primes $r\le63541$, define
\[
 \widehat D_r=\begin{cases}
 D_r,&r^{T_r}<2^{100},\\
 1,&r^{T_r}\ge2^{100}.
 \end{cases}
\]
For $29\le p\le63541$, we retain the distribution of $\prod_{r<p}\widehat D_r$ up to $H_0=100000$ and evaluate $L_p^{\mathrm{sc}}$ using~\eqref{eq:replacement-bound}. All moments are computed from the untruncated factors by~\eqref{eq:moments}.

Put
\[
 Y=\prod_{r\le63541}D_r,\qquad
 \widehat Y=\prod_{r\le63541}\widehat D_r,\qquad H(x)=\mathbb E(Y-x)_+.
\]
Let
\[
 M_1^-\le\mathbb EY\le M_1^+,\qquad
 \mathbb EY^2\le M_2^+,\qquad \mathbb EY^3\le M_3^+
\]
be moment bounds computed with directed rounding. On $t\le H_0$, let $P_t^+\ge\mathbb P(\widehat Y=t)$ and $P_t^-\le\mathbb P(Y=t)$.
The lower bounds $P_t^-$ are computed by discarding outcomes for which $r^{T_r}\ge2^{100}$ for at least one $r$, rather than replacing the corresponding factors by $1$.
Since $\widehat Y\le Y$, for $0\le x\le H_0$ we have
\[
 H(x)=\mathbb EY-x+\mathbb E(x-Y)_+
 \le\mathbb EY-x+\mathbb E(x-\widehat Y)_+
 \le M_1^+-x+\sum_{1\le t<x}(x-t)P_t^+.
\]
For $x>H_0$, we use
\[
 4(x-H_0)(y-x)_+\le(y-H_0)_+^2\qquad(y\ge0).
\]
The inequality is immediate for $y\le x$. For $y>x$, subtracting its left side from its right side gives $(y-2x+H_0)^2\ge0$. Moreover,
\[
 \mathbb E(Y-H_0)_+^2
 =\mathbb EY^2-2H_0\mathbb EY+H_0^2-\mathbb E(H_0-Y)_+^2.
\]
Using upper bounds for the positive terms and lower bounds for the subtracted terms, and applying~\eqref{eq:cubic}, we obtain
\[
 H(x)\le\min\left\{
 \frac{M_2^+-2H_0M_1^-+H_0^2-\sum_{1\le t<H_0}(H_0-t)^2P_t^-}{4(x-H_0)},
 \frac{4M_3^+}{27x^2}\right\}\qquad(x>H_0).
\]

For $63541<p\le10^7$, put
\[
 V_p=\prod_{63541<r<p}D_r.
\]
For $0\le i<65$ and $j\in\{0,1\}$, let $E_{i,j}(p)$ be the event that all $T_r$ with $63541<r<p$ are at most $2$, exactly $i$ equal $1$, and exactly $j$ equal $2$. On this event $V_p=2^i3^j$. Define
\[
 u_{i,j}(p)=\mathbb E[V_p\mathbf1_{E_{i,j}(p)}]
           =2^i3^j\mathbb P(E_{i,j}(p)).
\]
The events are disjoint, and $Y$ is independent of $V_p$. Hence, for $x\ge0$,
\[
\begin{aligned}
 \mathbb E(YV_p-x)_+
 &\le\sum_{i,j}u_{i,j}(p)H\left(\frac{x}{2^i3^j}\right)
       +\mathbb EY\left(\mathbb EV_p-\sum_{i,j}u_{i,j}(p)\right),\\
 L_p^{\mathrm{sc}}&=\frac{c_p}{p-1}
             \mathbb E\bigl(YV_p-\delta_p(p-1)\bigr)_+.
\end{aligned}
\]
We use the displayed bounds for $H$ when $i<16$ and $H\le\mathbb EY$ when $16\le i<65$. The subtracted sum in the first-moment term still includes every $i$ with $0\le i<65$. Positive terms use upper bounds and the subtracted first-moment sum uses lower bounds.

Initially $u_{0,0}=1$ and all other $u_{i,j}$ vanish. At prime $p$, we evaluate $L_p^{\mathrm{sc}}$ before including $T_p$, and then update the weights $u_{i,j}$ via
\[
\begin{aligned}
u'_{i,0}
 &=\left(1-\frac{c_p}{p}\right)u_{i,0}
   +\frac{2c_p(p-1)}{p^2}u_{i-1,0},\\
u'_{i,1}
 &=\left(1-\frac{c_p}{p}\right)u_{i,1}
   +\frac{2c_p(p-1)}{p^2}u_{i-1,1}
   +\frac{3c_p(p-1)}{p^3}u_{i,0},
\end{aligned}
\]
where terms with negative indices are zero. The full moments $\mathbb EV_p$ and $\mathbb E[(YV_p)^3]$ are computed separately without truncation.

There are $J=664579$ primes up to $10^7$. Beyond this point the third-moment multiplier is at most $1+14.001/(p-1)$, and $p_j-1>15j$ by~\eqref{eq:prime-fifteen}. Let $\mathcal M_3$ be the upper bound obtained by upward rounding for the untruncated third moment after $p_J$. Applying Lemma~\ref{lem:index-tail} with $N=J$, $C=32/27$, $a=14.001$, and $c=15$ gives
\[
 \sum_{p>10^7}L_p^{\mathrm{sc}}
 \le\frac{32000\mathcal M_3}{27\cdot225\cdot15999\,J^2}.
\]
At scale $2^{50}$, the resulting scaled upper bounds for $\sum L_p^{\mathrm{sc}}$ are
\begin{center}
\begin{tabular}{lr}
\toprule
Range & Upper bound for $2^{50}\sum L_p^{\mathrm{sc}}$\\
\midrule
$29\le p\le63541$ & $2400891806166066$\\
$63541<p\le10^7$ & $5562926093570$\\
$p>10^7$ & $125132341126$\\
\midrule
Total & $2406579864600762$\\
\bottomrule
\end{tabular}
\end{center}

Combining \eqref{eq:initial-union-bound}, \eqref{eq:small-divisor-reduction-bound}, \eqref{eq:early-reduction-bound}, \eqref{eq:late-reduction-bound}, and the scalar upper bound through \eqref{eq:small-divisor-reduction} and \eqref{eq:accumulated-reductions}, we find that the sum of the mass of the initial union and the contributions of all subsequent stages is at most
\[
 U_{\mathrm{tot}}^+=I_0+\frac{2406579864600762}{2^{50}}
       -\frac{1230344533175833263}{10^{18}}
       -\frac{8837268985393504}{10^{18}}
       -\frac{5740271089542}{2^{50}}.
\]
Clearing denominators gives
\begin{equation}\label{eq:general-final-bound}
 U_{\mathrm{tot}}^+<\frac{997}{1000}<1.
\end{equation}
Corollary~\ref{cor:criterion} rules out a covering with $m\ge14000$, proving Theorem~\ref{thm:general}.

\subsection{Squarefree moduli}\label{sec:squarefree}

For squarefree moduli, only the indicators $\mathbf1_{\{T_{p_j}\ge1\}}$ enter the comparison. These are independent Bernoulli variables.
Suppose for contradiction that every modulus exceeds $9$. The only squarefree moduli supported on $\{2,3\}$ are $2,3,6$, so no classes are introduced at the first two primes and we take $\delta_1=\delta_2=0$.

Let $I_j$ be independent Bernoulli variables with
\[
 \mathbb P(I_j=1)=a_j=\frac1{p_j(1-\delta_j)}.
\]
At stage $k$, the comparison integer $R=\prod_{j<k}p_j^{I_j}$ has divisor count $\tau(R)=2^{S_{k-1}}$, where $S_{k-1}=I_1+\cdots+I_{k-1}$.
Each subset of the preceding primes gives at most one modulus. For $p_k\in\{5,7\}$, the empty subset is excluded because its modulus is below $10$. For $p_k\ge11$, every subset is admissible. Theorem~\ref{thm:scalar-comparison} and~\eqref{eq:stage-measure} give
\begin{equation}\label{eq:squarefree-contribution}
 L_k\le a_k\mathbb E\left(2^{S_{k-1}}-p_k\delta_k-\mathbf1_{\{3\le k\le4\}}\right)_+,
\end{equation}
where $L_k$ is the mass introduced at $p_k$.

Put
\[
 u_{k,s}=\mathbb P(S_{k-1}=s),\qquad
 M_k=\mathbb E2^{S_{k-1}},\qquad R_k=\mathbb E8^{S_{k-1}}.
\]
Independence gives
\begin{equation}\label{eq:squarefree-recursion}
\begin{aligned}
 u_{k+1,0}&=(1-a_k)u_{k,0},\qquad
 u_{k+1,s}=(1-a_k)u_{k,s}+a_ku_{k,s-1}\quad(s\ge1),\\
 M_{k+1}&=M_k(1+a_k),\qquad
 R_{k+1}=R_k(1+7a_k),
\end{aligned}
\end{equation}
with $u_{1,0}=M_1=R_1=1$ and $u_{1,s}=0$ for $s>0$.
For $A\ge0$, \eqref{eq:positive-part} gives
\begin{equation}\label{eq:squarefree-positive}
 \mathbb E(2^{S_{k-1}}-A)_+
 =M_k-A+\sum_{2^s<A}(A-2^s)u_{k,s}.
\end{equation}
For $N=10^6$, the parameter choice in Table~\ref{tab:squarefree-parameters} satisfies
\[
 p_k\delta_k+\mathbf1_{\{3\le k\le4\}}<8\cdot10^6<2^{23}
       \qquad(3\le k\le N).
\]
Thus values $0\le s\le40$ suffice. We evaluate $1-a_k$ exactly and round all updates in~\eqref{eq:squarefree-recursion} upward at scale $10^{45}$.

\begin{proposition}\label{prop:squarefree-numerics}
With the parameters in Table~\ref{tab:squarefree-parameters},
\begin{equation}\label{eq:squarefree-finite}
 \sum_{k=3}^N L_k\le\SquarefreeFiniteUpper
 \qquad\text{and}\qquad R_{N+1}\le\SquarefreeThirdMomentUpper.
\end{equation}
\end{proposition}

\begin{proof}
For $k=1,\ldots,N$, we compute the next distribution and moments from~\eqref{eq:squarefree-recursion}. For $k\ge3$, before this update we substitute~\eqref{eq:squarefree-positive}, with $A=p_k\delta_k+\mathbf1_{\{3\le k\le4\}}$, into~\eqref{eq:squarefree-contribution}. At scale $10^{45}$, upward rounding gives~\eqref{eq:squarefree-finite}. The first $500$ steps were also checked independently in exact rational arithmetic.
\end{proof}

For $k>N$, take $\delta_k=1/2$. Then \eqref{eq:cubic} gives $L_k\le32R_k/(27p_k^3)$, while $R_{k+1}=R_k(1+14/p_k)$ and $p_k>15k$ by~\eqref{eq:prime-fifteen}. Lemma~\ref{lem:index-tail}, with $C=32/27$, $a=14$, and $c=15$, yields
\[
 \sum_{k>N}L_k\le\frac{2R_{N+1}}{6075N^2}\le\SquarefreeTailUpper.
\]
Combining the finite and tail numerators at scale $10^{45}$ gives
\[
 \sum_{k\ge3}L_k\le\frac{\SquarefreeTotalNumerator}{10^{45}}
                  \le\SquarefreeTotalUpper<1.
\]
Corollary~\ref{cor:criterion} rules out a covering in which all squarefree moduli exceed $9$. Since neither $8$ nor $9$ is squarefree, some modulus is at most $7$, proving Theorem~\ref{thm:squarefree}.

\appendix
\section{Computational details}\label{app:computations}

\subsection{Computational verification and rounding}\label{sec:verification-rounding}

The finite recursions in Sections~\ref{sec:distortion-scalar-comparison}--\ref{sec:bounds-proof} are evaluated using C++17 and Python programs, with 128-bit and GMP/Boost multiprecision integer arithmetic at the stated scales $W$.

Upper bounds for the stage contributions are rounded upward, whereas lower bounds for the subtracted reductions $G_p$, $\mathcal R_{[29,233]}$, and $\mathcal R_{>233}$ are rounded downward. Moment bounds are rounded according to their role in the estimate. When computing an upper bound, positive terms are bounded from above and subtracted terms from below. The directions are reversed when computing a lower bound. Differences of positive parts are put over a common denominator before rounding. At each stage, the contribution is evaluated before the distribution for the next prime is updated. The final rational non-covering inequalities are verified by exact integer multiplication.

Independent Python implementations were used to recompute the complete integer recurrences for Theorem~\ref{thm:odd-lcm} and Theorem~\ref{thm:squarefree}. For Theorem~\ref{thm:odd}, the three cases $m\ge15$, $m=11$, and $m=13$ were evaluated by C++ programs. For Theorem~\ref{thm:general}, the programs generate all $23$-smooth integers up to $10^{12}$, of which there are $892644$, and the $512$ Shearer support polynomials, evaluate the scalar bound and the small-divisor and comparison-prime reductions, and assemble~\eqref{eq:general-final-bound} using exact arithmetic. Independent checks were also carried out for some intermediate calculations and for the final arithmetic.

\subsection{Fixed rational parameters}\label{app:parameters}

\emph{General case.}
The parameter $\delta_p$ is $0$ for every prime up to $23$. For the $6360$ primes from $29$ to $63541$, the full parameter table records $10^6\delta_p$ as an integer. The first twelve values are shown in the following table.
\begin{center}
\begin{tabular}{rr@{\qquad}rr}
\toprule
$p$ & $10^6\delta_p$ & $p$ & $10^6\delta_p$\\\midrule
$29$ & $30000$ & $53$ & $99077$\\
$31$ & $40000$ & $59$ & $106093$\\
$37$ & $36929$ & $61$ & $118857$\\
$41$ & $55562$ & $67$ & $124385$\\
$43$ & $74789$ & $71$ & $131391$\\
$47$ & $72667$ & $73$ & $141938$\\
\bottomrule
\end{tabular}
\end{center}
For larger primes, we set $\delta_p=41323/100000$ for $63541<p\le10^7$ and $\delta_p=1/2$ for $p>10^7$.
We also take $\vartheta_p=1/c_p$ for $29\le p\le63541$.

\medskip
\noindent\emph{Squarefree case.}
Table~\ref{tab:squarefree-parameters} gives the complete choice of parameters for the primes $p_k$.
\begin{table}[!htbp]
\centering
\caption{Parameters for squarefree moduli}\label{tab:squarefree-parameters}
\begin{tabular}{cr@{\qquad\qquad}cr}
\toprule
Prime index $k$ & $\delta_k$ & Prime index $k$ & $\delta_k$\\\midrule
$1$--$2$ & $0$ & $97$ & $1/2$\\
$3$--$4$ & $1/p_k$ & $98$--$166$ & $256/p_k$\\
$5$--$8$ & $4/p_k$ & $167$--$172$ & $1/2$\\
$9$--$12$ & $8/p_k$ & $173$--$292$ & $512/p_k$\\
$13$--$20$ & $16/p_k$ & $293$--$309$ & $1/2$\\
$21$--$33$ & $32/p_k$ & $310$--$500$ & $1024/p_k$\\
$34$--$56$ & $64/p_k$ & $k\ge501$ & $1/2$\\
$57$--$96$ & $128/p_k$ &  & \\
\bottomrule
\end{tabular}
\end{table}

\noindent\emph{Odd moduli and the least common multiple.}
Table~\ref{tab:odd-distortion-parameters} specifies the integer $A_p$ in $\delta_p=A_p/(p-1)\le1/2$ for $5\le p\le7993$. These parameters are used from $p=5$ onward in the proof of $9\mid Q$ and after the initial prime steps in Section~\ref{sec:odd}. For larger primes, we set $\delta_p=71/200$ for $7993<p\le10^5$ and $\delta_p=1/2$ for $p>10^5$. The proof of $5\mid Q$ uses $\delta_p=3/10$ for all $p\ge7$, with $\gamma_3=2$.
\begin{table}[htbp]
\centering\small
\caption{Integer parameters $A_p$ for the distortion estimates}\label{tab:odd-distortion-parameters}
\begin{tabular}{cr@{\qquad}cr@{\qquad}cr}
\toprule
Prime range & $A_p$ & Prime range & $A_p$ & Prime range & $A_p$\\\midrule
$5$ & $1$ & $307$ & $120$ & $2161$ & $960$\\
$7\le p\le11$ & $2$ & $311\le p\le419$ & $128$ & $2179\le p\le2843$ & $1024$\\
$13\le p\le19$ & $4$ & $421\le p\le449$ & $144$ & $2851\le p\le3163$ & $1152$\\
$23$ & $6$ & $457$ & $160$ & $3167\le p\le3209$ & $1280$\\
$29\le p\le31$ & $8$ & $461\le p\le577$ & $192$ & $3217$ & $1296$\\
$37$ & $9$ & $587$ & $216$ & $3221\le p\le4099$ & $1536$\\
$41\le p\le43$ & $12$ & $593\le p\le787$ & $256$ & $4111\le p\le4153$ & $1728$\\
$47\le p\le61$ & $16$ & $797\le p\le859$ & $288$ & $4157$ & $1792$\\
$67$ & $18$ & $863$ & $320$ & $4159\le p\le4177$ & $1920$\\
$71\le p\le83$ & $24$ & $877\le p\le1103$ & $384$ & $4201\le p\le5443$ & $2048$\\
$89\le p\le113$ & $32$ & $1109\le p\le1117$ & $432$ & $5449\le p\le6053$ & $2304$\\
$127$ & $40$ & $1123$ & $480$ & $6067\le p\le6151$ & $2560$\\
$131\le p\le157$ & $48$ & $1129\le p\le1493$ & $512$ & $6163$ & $2592$\\
$163\le p\le223$ & $64$ & $1499\le p\le1637$ & $576$ & $6173$ & $2880$\\
$227\le p\le233$ & $72$ & $1657\le p\le1667$ & $640$ & $6197\le p\le7879$ & $3072$\\
$239$ & $80$ & $1669\le p\le2131$ & $768$ & $7883\le p\le7993$ & $3456$\\
$241\le p\le293$ & $96$ & $2137\le p\le2153$ & $864$ &  & \\
\bottomrule
\end{tabular}
\end{table}

\subsection{Rational bounds for the exponential estimates}\label{app:exponential-checks}

For $0<x<n+2$, the Taylor series and a geometric series bound for the remainder give
\begin{equation*}
e^x\le B_n(x):=\sum_{k=0}^{n}\frac{x^k}{k!}
+\frac{x^{n+1}}{(n+1)!}\frac1{1-x/(n+2)},
\end{equation*}
since the ratio of successive terms after degree $n+1$ is at most $x/(n+2)<1$. Evaluating $B_n$ at $n=64$ and $n=24$ by exact rational arithmetic gives
\begin{equation*}
B_{64}\!\left(\frac{6703}{500}\right)<663976<664579,\qquad
B_{24}\!\left(\frac{519}{200}\right)<\frac{67}{5}<\frac{6703}{500},
\end{equation*}
which establishes~\eqref{eq:exponential-bounds} since $6703/500=13.406$ and $519/200=2.595$.

\section*{Data and code availability}
The code and data used for the computations in this paper are available from the authors upon request.

\section*{AI Use Disclosure}
The research direction and initial mathematical framework were developed by the authors.
Generative AI was used interactively in the subsequent mathematical exploration, including
the development, testing, and refinement of parts of the construction and proof. The authors
independently verified the final arguments and take full responsibility for the paper.

\section*{Acknowledgments}
The authors would like to thank Biao Wang for helpful discussions.
Huixi Li is supported by the National Natural Science Foundation of China (Grant No. 12561001).
Shaoyun Yi is supported by the National Natural Science Foundation of China (No. 12471187) and the Fundamental Research Funds for the Central Universities (No. 20720260221).

\begingroup
\small
\setlength{\bibsep}{2pt}
\needspace{8\baselineskip}

\endgroup


\begin{thebibliography}{10}

\bibitem{Azlin2011thesis}
Michael~Wayne Azlin.
\newblock Covering systems of polynomial rings over finite fields.
\newblock Master's thesis, University of Mississippi, 2011.

\bibitem{Balister2024}
Paul Balister.
\newblock Erd\H{o}s covering systems.
\newblock In {\em Surveys in combinatorics 2024}, volume 493 of {\em London
  Math. Soc. Lecture Note Ser.}, pages 31--54. Cambridge Univ. Press,
  Cambridge, 2024.

\bibitem{BBMST2020}
Paul Balister, B\'{e}la Bollob\'{a}s, Robert Morris, Julian Sahasrabudhe, and
  Marius Tiba.
\newblock Erd\H{o}s covering systems.
\newblock {\em Acta Math. Hungar.}, 161(2):540--549, 2020.

\bibitem{BBMST2021squarefree}
Paul Balister, B\'ela Bollob\'as, Robert Morris, Julian Sahasrabudhe, and
  Marius Tiba.
\newblock The {E}rd\H{o}s--{S}elfridge problem with square-free moduli.
\newblock {\em Algebra Number Theory}, 15(3):609--626, 2021.

\bibitem{BBMST2022}
Paul Balister, B\'{e}la Bollob\'{a}s, Robert Morris, Julian Sahasrabudhe, and
  Marius Tiba.
\newblock On the {E}rd{\H{o}}s covering problem: the density of the uncovered
  set.
\newblock {\em Invent. Math.}, 228(1):377--414, 2022.

\bibitem{BFF1986}
Marc A. Berger, Alexander Felzenbaum, and Aviezri S. Fraenkel.
\newblock Necessary condition for the existence of an incongruent covering system with odd moduli.
\newblock {\em Acta Arith.}, 45(4):375--379, 1986.

\bibitem{BFF1987}
Marc A. Berger, Alexander Felzenbaum, and Aviezri S. Fraenkel.
\newblock Necessary condition for the existence of an incongruent covering system with odd moduli. {II}.
\newblock {\em Acta Arith.}, 48(1):73--79, 1987.

\bibitem{BispelsEtAl2026}
Chris Bispels, Matthew Cohen, Joshua Harrington, Joshua Lowrance, Kaelyn
  Pontes, Leif Schaumann, and Tony W.~H. Wong.
\newblock A further investigation on covering systems with odd moduli.
\newblock {\em Discrete Math.}, 349:115013, 2026.

\bibitem{BloomErdosProblems}
T.~F. Bloom.
\newblock Erd\H{o}s Problems, Problems 2 and 7.
\newblock \url{https://www.erdosproblems.com/2} and
  \url{https://www.erdosproblems.com/7}.
\newblock Accessed 10 November 2026.

\bibitem{Choi1971}
S.~L.~G. Choi.
\newblock Covering the set of integers by congruence classes of distinct
  moduli.
\newblock {\em Math. Comp.}, 25:885--895, 1971.

\bibitem{Churchhouse1968}
R.~F. Churchhouse.
\newblock Covering sets and systems of congruences.
\newblock In R.~F. Churchhouse and J.-C. Herz, editors, {\em Computers in
  Mathematical Research}, pages 20--36. North-Holland, Amsterdam, 1968.

\bibitem{CFT2025}
Maria Cummings, Michael Filaseta, and Ognian Trifonov.
\newblock An upper bound for the minimum modulus in a covering system with
  squarefree moduli.
\newblock {\em Acta Math. Hungar.}, 175(1):1--25, 2025.

\bibitem{Dai2026}
Xiangjun Dai.
\newblock Covering systems of polynomial rings over finite fields.
\newblock Master's thesis, Nankai University, 2026.
\newblock In Chinese.

\bibitem{Dusart2010}
Pierre Dusart.
\newblock Estimates of some functions over primes without {R.H.}, 2010.
\newblock arXiv:1002.0442.

\bibitem{Erdos1950}
Paul Erd{\H{o}}s.
\newblock On integers of the form {$2^k+p$} and some related problems.
\newblock {\em Summa Brasil. Math.}, 2:113--123, 1950.

\bibitem{Erdos1995}
Paul Erd\H{o}s.
\newblock Some of my favourite problems in number theory, combinatorics, and
  geometry.
\newblock {\em Resenhas}, 2(2):165--186, 1995.
\newblock Combinatorics Week (Portuguese) (S\~{a}o Paulo, 1994).

\bibitem{FFKPY2007}
Michael Filaseta, Kevin Ford, Sergei Konyagin, Carl Pomerance, and Gang Yu.
\newblock Sieving by large integers and covering systems of congruences.
\newblock {\em J. Amer. Math. Soc.}, 20(2):495--517, 2007.

\bibitem{FK2026}
Michael Filaseta and Alexandros Kalogirou.
\newblock Covering systems with the sum of the reciprocals of the moduli close to $1$.
\newblock {\em Trans. Amer. Math. Soc.}, 379:4945--4972, 2026.

\bibitem{FKK2026}
Michael Filaseta, Alexandros Kalogirou, and Jonah Klein.
\newblock Results for distinct covering systems with a prime minimum modulus $q\ge11$.
\newblock 2026.
\newblock arXiv:2610.34541v1.

\bibitem{Gibson2009}
Donald~Jason Gibson.
\newblock A covering system with least modulus 25.
\newblock {\em Math. Comp.}, 78(266):1127--1146, 2009.

\bibitem{GuoSun2005}
Song Guo and Zhi-Wei Sun.
\newblock On odd covering systems with distinct moduli.
\newblock {\em Adv. in Appl. Math.}, 35(2):182--187, 2005.

\bibitem{Gyarmati2026}
Katalin Gyarmati.
\newblock The case of the second smallest modulus: or why we cannot be arbitrarily careless?, 2026.
\newblock arXiv:2606.05251v1.

\bibitem{HarringtonSunWong2022}
Joshua Harrington, Yewen Sun, and Tony W.~H. Wong.
\newblock Covering systems with odd moduli.
\newblock {\em Discrete Math.}, 345(8):112936, 2022.

\bibitem{Hough2015}
Bob Hough.
\newblock Solution of the minimum modulus problem for covering systems.
\newblock {\em Ann. of Math. (2)}, 181(1):361--382, 2015.

\bibitem{HoughNielsen2019}
Robert~D. Hough and Pace~P. Nielsen.
\newblock Covering systems with restricted divisibility.
\newblock {\em Duke Math. J.}, 168(17):3261--3295, 2019.

\bibitem{KKL2024}
Jonah Klein, Dimitris Koukoulopoulos, and Simon Lemieux.
\newblock On the {$j$}-th smallest modulus of a covering system with distinct moduli.
\newblock {\em Int. J. Number Theory}, 20(2):471--479, 2024.

\bibitem{Krukenberg1971}
Claire~Emil Krukenberg.
\newblock {\em Covering sets of the integers}.
\newblock PhD thesis, University of Illinois at Urbana-Champaign, 1971.

\bibitem{LWWY2024}
Huixi Li, Biao Wang, Chunlin Wang, and Shaoyun Yi.
\newblock On covering systems of polynomial rings over finite fields.
\newblock {\em Proc. Amer. Math. Soc.}, 152(9):3731--3742, 2024.

\bibitem{LWWY2025}
Huixi Li, Biao Wang, Chunlin Wang, and Shaoyun Yi.
\newblock On {E}rd\H{o}s covering systems in global function fields.
\newblock {\em J. Number Theory}, 266:269--280, 2025.

\bibitem{LiWangYi2023pre}
Huixi Li, Biao Wang, and Shaoyun Yi.
\newblock On the minimum modulus problem in number fields, 2023.
\newblock arXiv:2302.05946, v2 (2026).

\bibitem{Morikawa1981examples}
Ryozo Morikawa.
\newblock Some examples of covering sets.
\newblock {\em Bull. Fac. Liberal Arts Nagasaki Univ. Natur. Sci.}, 21(2):1--4,
  1981.

\bibitem{Nie2026}
Junzhi Nie.
\newblock Modular transfer graphs and the construction of covering systems over
  {$\mathbb F_q[x]$}.
\newblock Master's thesis, Nankai University, 2026.
\newblock In Chinese.

\bibitem{Nielsen2009}
Pace~P. Nielsen.
\newblock A covering system whose smallest modulus is 40.
\newblock {\em J. Number Theory}, 129(3):640--666, 2009.

\bibitem{Owens2014}
Tyler Owens.
\newblock A covering system with minimum modulus 42.
\newblock Master's thesis, Brigham Young University, 2014.

\bibitem{Rosser1939}
Barkley Rosser.
\newblock The {$n$}th prime is greater than {$n\log n$}.
\newblock {\em Proc. London Math. Soc. (2)}, 45(1):21--44, 1939.

\bibitem{Schinzel1967}
A.~Schinzel.
\newblock Reducibility of polynomials and covering systems of congruences.
\newblock {\em Acta Arith.}, 13:91--101, 1967.

\bibitem{Shearer1985}
J.~B. Shearer.
\newblock On a problem of {S}pencer.
\newblock {\em Combinatorica}, 5(3):241--245, 1985.

\bibitem{SimpsonZeilberger1991}
R.~J. Simpson and Doron Zeilberger.
\newblock Necessary conditions for distinct covering systems with square-free
  moduli.
\newblock {\em Acta Arith.}, 59(1):59--70, 1991.

\bibitem{Sun1991}
Zhi-Wei Sun.
\newblock On covering systems with distinct moduli.
\newblock {\em J. Yangzhou Teachers College Nat. Sci. Ed.}, 11(3):21--27, 1991.
\newblock In Chinese.

\bibitem{Swift1954}
J.~D. Swift.
\newblock Sets of covering congruences.
\newblock {\em Bull. Amer. Math. Soc.}, 60(4):390, 1954.
\newblock Abstract 561.

\bibitem{Tchen1980}
Andr\'e{}~H. Tchen.
\newblock Inequalities for distributions with given marginals.
\newblock {\em Ann. Probab.}, 8(4):814--827, 1980.

\bibitem{BiaoWang2026}
Biao Wang.
\newblock Bounds for {Erd\H{o}s} covering systems in global function fields.
\newblock {\em Int. J. Number Theory}, 22(6):1327--1336, 2026.

\bibitem{RongyinWang2026}
Rongyin Wang.
\newblock On an {Erd\H{o}s}-type conjecture on {$\mathbb F_q[x]$}.
\newblock {\em Finite Fields Appl.}, 109:102720, 2026.

\bibitem{RongyinWang2026asymptotic}
Rongyin Wang.
\newblock An asymptotic bound for non-covering congruence systems over
  {$\mathbb F_q[x]$}, 2026.
\newblock arXiv:2607.27538, v2.

\end{thebibliography}
\end{document}